\documentclass[11pt]{article}

\usepackage{epsfig}
\usepackage{wrapfig}
\input epsf
\usepackage{here}
\usepackage{latexsym}

\usepackage{graphicx}
\usepackage[dvips]{color}
\usepackage{latexsym}
\usepackage{amsmath}
\usepackage{graphics}
\usepackage{amsthm}
\usepackage{amssymb}

\usepackage
{subcaption}

  \makeatletter
  \@addtoreset{footnote}{page}
  \makeatother

\makeatletter

\advance\headheight1.2pt

\newcommand\shavedisplay{%
	\abovedisplayskip=4pt
	\belowdisplayskip=6pt
	\jot=0pt
}

\newtheorem{theorem}{Theorem}[section]

\newtheorem{proposition}{Proposition}[section]
\newtheorem{lemma}{Lemma}[section]

\newtheorem{remark}{Remark}[section]
\newtheorem{definition}{Definition}[section]

\newtheorem{conjecture}{Conjecture}[section]

\newtheorem{o-problem}{Open question}[section]

\newtheorem{fact}{Fact}[section]

\begin{document}
\shavedisplay
\pagestyle{plain}

  \title{Non-uniqueness of closed embedded non-smooth hypersurfaces with constant anisotropic mean curvature}

\author{Yoshiki Jikumaru and Miyuki Koiso\footnote{
This work was partially supported by JSPS KAKENHI Grant Number JP18H04487. 
}}

\date{}


  \maketitle

\begin{abstract}
An anisotropic surface energy is the integral of an energy density that depends on the normal at each point over the considered surface, and it is a generalization of surface area. The minimizer of such an energy among all closed surfaces enclosing the same volume is unique and it is (up to rescaling) so-called the Wulff shape. We prove that, unlike the isotropic case,  
there exists an 
anisotropic energy density function such that there exist closed embedded equilibrium surfaces with genus zero in ${\mathbb R}^3$ each of which is not (any homothety and translation of) the Wulff shape. 
We also give nontrivial self-similar shrinking solutions of anisotropic mean curvature flow. 
These results are generalized to hypersurfaces in ${\mathbb R}^{n+1}$. 
\end{abstract}


\section{Introduction}\label{intro}

An anisotropic surface energy
is the integral of an energy density that depends on the surface normal, which was introduced by J. W. Gibbs (1839-1903) in order to model the shape of small crystals (\cite{W},\cite{Wu}). 
The energy density function is given as 
a positive continuous function $\gamma:S^n \rightarrow {\mathbb R}_{> 0}$ on the unit sphere $S^n=\{\nu \in {\mathbb R}^{n+1} \; | \; \|\nu\|=1\}$ in the $(n+1)$-dimensional Euclidean space ${\mathbb R}^{n+1}$. 
In this paper, for a technical reason, we study the case where $\gamma$ is of $C^2$. 
When we study variational problems of such an energy, it is natural to consider not only smooth surfaces but also surfaces with edges. 
And so let  
$X:M\to {\mathbb R}^{n+1}$ be 
a piecewise-$C^2$ weak immersion (the definition of piecewise-$C^2$ weak immersion will be given in \S \ref{pw}) of 
an $n$-dimensional oriented connected compact $C^\infty$ manifold $M$ into ${\mathbb R}^{n+1}$. 
The unit normal vector field $\nu$ along $X$ is defined on $M$ except a set with measure zero. The anisotropic energy ${\cal F}_\gamma(X)$ of $X$ is defined as  
$
{\cal F}_\gamma(X)=\int_{M} \gamma(\nu)\:dA$, 
where $dA$ is the $n$-dimensional volume form of $M$ induced by $X$. 
In the special case where $\gamma\equiv1$, ${\mathcal F}_\gamma(X)$ is the usual $n$-dimensional volume of $X$ that is the volume of $M$ with metric induced by $X$. 

For any positive number $V>0$, among all closed hypersurfaces in ${\mathbb R}^{n+1}$ enclosing the same $(n+1)$-dimensional volume $V$, there exists a unique (up to translation in ${\mathbb R}^{n+1}$) minimizer of ${\mathcal F}_\gamma$ (\cite{T}). 
Here a closed hypersurface means that the boundary (having tangent space almost everywhere) of a set of positive Lebesgue measure. The minimizer for 
$V=(n+1)^{-1}\int_{S^n} \gamma(\nu) \;dS^n$ is called the Wulff shape (for $\gamma$) (see \S \ref{con} for the standard definition of the Wulff shape), and we denote it by $W_\gamma$. $W_1$ is the unit sphere $S^n$. 
In general, $W_\gamma$ is convex but not necessarily smooth. 

A piecewise-$C^2$ weak immersion $X:M\to {\mathbb R}^{n+1}$ is a critical point of ${\mathcal F}_\gamma$ for all variations that preserve the enclosed $(n+1)$-dimensional volume (we will call such a variation a volume-preserving variation) if and only if the anisotropic mean curvature of $X$ is constant and $X$ has a certain condition on its singular points (cf. \S \ref{CHM}). 
Here the anisotropic mean curvature $\Lambda$ of $X$ is defined at each regular point of $X$ as
$
\Lambda:=(1/n) \large(-{\rm div}_M D\gamma+nH\gamma\large)$, 
where $D\gamma$ is the gradient of $\gamma$ on $S^n$ and $H$ is the mean curvature of $X$ (cf. \cite{R}, \cite{KP2005}).  
We call such $X$ a CAMC (constant anisotropic mean curvature) hypersurface (Definition \ref{defR}). When $\gamma\equiv 1$, $\Lambda$ coincides with $H$ and CAMC hypersurfaces are CMC (constant mean curvature) hypersurfaces (of $C^\omega$ class).

A natural question is 
`is any closed CAMC hypersurface $X$ homothetic to the Wulff shape?' 
The answer is not affirmative even in the case where $\gamma\equiv 1$ (\cite{HW1986}, \cite{Kap1990}, \cite{Kap1995}). 
However, it is expected that, if $X$ satisfies one of the following conditions (I)-(III), the image of $X$ is a homothety of the Wulff shape. 


(I) $X$ is an embedding, that is, $X$ is an injective mapping.

(II) $X$ is stable.

(III) $n=2$ and the genus of $M$ is $0$, that is, $M$ is homeomorphic to $S^2$. 


\noindent 
Here a CAMC hypersurface is said to be stable if the second variation of the energy ${\mathcal F}_\gamma$ for any volume-preserving variation is nonnegative. 
If we assume that the Wulff shape $W_\gamma$ is a smooth strictly convex hypersurface, any closed CAMC hypersurface $X$ is also smooth and the above expectation is correct, which was proved by \cite{A}, \cite{BD1984}, \cite{H} for $\gamma\equiv 1$, and by \cite{HLMG2009}, \cite{P1998}, \cite{KP2010} and \cite{Ando2012} for general $\gamma$.

If $\gamma$ has less regularity or less convexity, the Wulff shape and CAMC hypersurfaces can have ``edges''.  
The second author \cite{K2018} studies the case where $\gamma:S^n\to {\mathbb R}_{>0}$ is a $C^2$ function and proved that, if $\gamma$ is convex (see \S \ref{con} for the definition of convexity), then 
any closed stable CAMC hypersurface for $\gamma$ 
whose $r$-th anisotropic mean 
curvature for $\gamma$ (Definition \ref{curv}) is integrable for $r=1, \cdots, n$ is (up to translation and homothety) the Wulff shape. 
However, the situation is not the same for non-convex $\gamma$. 
Actually in this paper we prove the following non-uniqueness results. 

\begin{theorem}\label{nonuni1}
There exists a $C^\infty$ function $\gamma:S^n\to {\mathbb R}_{>0}$ which is not a convex integrand such that there exist closed embedded CAMC hypersurfaces in ${\mathbb R}^{n+1}$ for $\gamma$ each of which is not (any homothety or translation of) the Wulff shape. 
\end{theorem}

\begin{theorem}\label{nonuni2}
There exists a $C^\infty$ function $\gamma:S^2\to {\mathbb R}_{>0}$ which is not a convex integrand such that there exist closed embedded CAMC surfaces in ${\mathbb R}^3$ with genus zero for $\gamma$ each of which is not (any homothety or translation of) the Wulff shape. 
\end{theorem}

These results are proved by giving  examples (\S \ref{ex1}, \ref{ex2}, \ref{pf}). 

The same examples will be applied to the anisotropic mean curvature flow. 
In order to give the precise statement, we recall that 
the Cahn-Hoffman map $\xi_\gamma$ for $\gamma$ is the mapping $\xi_\gamma:S^n \to {\mathbb R}^{n+1}$ defined as
$
\xi_\gamma(\nu)=D\gamma|_{\nu}+\gamma(\nu)\nu$, here the tangent space $T_\nu(S^n)$ of $S^n$ at $\nu \in S^n$ is naturally identified with the $n$-dimensional linear subspace of ${\mathbb R}^{n+1}$ (cf. \cite{K2018}). 
Let $X_t:M \to {\mathbb R}^{n+1}$ be one-parameter family of piecewise-$C^2$ weak immersions with anisotropic mean curvature $\Lambda_t$. Assume that the Cahn-Hoffman field $\tilde{\xi}_t$ (which is an anisotropic generalization of the unit normal vector field. see \S \ref{CHM} for definition) along $X_t$ is defined on $M$. 
If $X_t$ satisfies 
$\partial X_t/\partial t = \Lambda_t \tilde{\xi}_t$, 
it is called an anisotropic mean curvature flow, which diminishes the anisotropic energy if $\Lambda_t \not\equiv0$ (\S \ref{app}). 
By a simple observation we show the following (\S \ref{app}).

\begin{theorem}\label{ASS}
Let $c$ be a positive constant. Set 
$$
X_t:=\sqrt{2(c-t)}\;\xi_\gamma, \quad t\le c.
$$
Then $X_t$ is a self-similar shrinking solution of the anisotropic mean curvature flow for $\gamma$, that is 

(i) $\partial X_t/\partial t = \Lambda_t \tilde{\xi}_t$, and 

(ii) $X_t$ is homothetic to $\xi_\gamma$ and it shrinks as $t$ increases. 
\end{theorem}

By using this result and by giving examples, we prove the following result (\S \ref{ex1}, \ref{ex2}, \ref{app}).

\begin{theorem}\label{AMCF}
There exists a $C^\infty$ function $\gamma:S^n\to {\mathbb R}_{>0}$ which is not a convex integrand such that there exist closed embedded self-similar shrinking solutions in ${\mathbb R}^{n+1}$ for $\gamma$ each of which is homeomorphic to $S^n$ and is not (any homothety or translation of) the Wulff shape. 
\end{theorem}
In contrast with this result, the round sphere is the only closed embedded self-similar shrinking solution of mean curvature flow in ${\mathbb R}^3$ with genus zero 
(\cite{SB2016}). 

Finally we give two conjetures about the uniqueness problems studied in Theorems \ref{nonuni1}, \ref{nonuni2}, and \ref{AMCF}. 

\begin{conjecture}\label{conj1}
Assume that $\gamma:S^n\to {\mathbb R}_{>0}$ is of $C^2$. 
Let 
$X:M\to {\mathbb R}^{n+1}$ be 
a closed CAMC hypersurface. 
We assume that the $r$-th anisotropic mean 
curvature of $X$ for $\gamma$ is  integrable for $r=1, \cdots, n$. Then, if $X$ satisfies at least one of the conditions (I), (II), (III) above, $X(M)$ is a subset of a homothety of the image $\xi_\gamma(S^n)$ of the Cahn-Hoffman map $\xi_\gamma$.
\end{conjecture}

\begin{conjecture}\label{conj2}
Assume that $\gamma:S^2\to {\mathbb R}_{>0}$ is of $C^2$. 
Then any closed embedded self-similar shrinking solution of the anisotropic mean curvature flow for $\gamma$ in ${\mathbb R}^{3}$ with genus zero is a subset of a homothety of $\xi_\gamma(S^2)$.
\end{conjecture}

This article is organized as follows. 
\S \ref{pre} is the section of preliminaries. 
In \S \ref{ex1} we give an example which proves Theorem \ref{nonuni1} for $n=1$ (Proposition \ref{prop}). In \S \ref{ex2} we give examples which will be used to prove Theorems \ref{nonuni1}, \ref{nonuni2}, and \ref{AMCF}. 
In \S \ref{pf} we give proofs of Theorems \ref{nonuni1}, \ref{nonuni2}. 
In \S \ref{app}, 
we prove Theorems \ref{ASS},  \ref{AMCF}.

\section{Preliminaries} \label{pre}

\subsection{Definitions of piecewise-$C^r$ weak immersion and its anisotropic energy} \label{pw}

First we recall the definition of a {\it piecewise-$C^r$ weak immersion}, ($r\in {\mathbb N}$), defined in \cite{K2018}. 
Let 
$M=\cup_{i=1}^k M_i$ be an $n$-dimensional oriented compact connected $C^\infty$ manifold, where each $M_i$ is an $n$-dimensional connected compact submanifold of $M$ with boundary, and $M_i\cap M_j = \partial M_i\cap\partial M_j$, ($i, j \in \{1, \cdots, k\}$, $i\ne j$). 
We call a map 
$X:M \rightarrow{\mathbb R}^{n+1}$ a piecewise-$C^r$ weak immersion (or a piecewise-$C^r$ weakly immersed hypersurface) if $X$  satisfies the following conditions (A1), (A2), and (A3) ($i \in \{1, \cdots, k\}$). 

(A1) $X$ is continuous, and each $X_i:=X|_{M_i} : M_i \to {\mathbb R}^{n+1}$ is of $C^r$.

(A2) The restriction $X|_{M_i^o}$ of $X$ to the interior $M_i^o$ of $M_i$ is a $C^r$-immersion. 

(A3) The unit normal vector field $\nu_i:{M_i^o} \to S^n$ along $X_i|_{M_i^o}$ can be extended to a $C^{r-1}$-mapping $\nu_i:M_i \to S^n$. 
Here, if $(u^1, \cdots, u^n)$ is a local coordinate system in $M_i$, $\{\nu_i, \partial/\partial u^1, \cdots,\partial/\partial u^n\}$ gives the canonical orientation in ${\mathbb R}^{n+1}$. 


\vskip0.3truecm
The image $X(M)$ of a piecewise-$C^r$ weak immersion $X:M \rightarrow{\mathbb R}^{n+1}$ is also called a piecewise-$C^r$ hypersurface. 
Denote by $S(X)$ the set of all singular points of $X$, here a singular point of $X$ is a point in $M$ at which $X$ is not an immersion. 


The anisotropic energy of a piecewise-$C^1$  weak immersion $X:M \rightarrow{\mathbb R}^{n+1}$ is defined as follows. 
Assume that $\gamma:S^n \rightarrow {\mathbb R}_{\ge 0}$ is a nonnegative continuous function. 
Let 
$\nu:M\setminus S(X) \to S^n$ be the unit normal vector field along $X|_{M\setminus S(X)}$. The anisotropic energy ${\cal F}_\gamma(X)$ of $X$ is defined as  
\begin{equation}\label{ae}
{\cal F}_\gamma(X):=\int_{M} \gamma(\nu)\:dA :=\sum_{i=1}^k\int_{M_i} \gamma(\nu_i)\:dA.
\end{equation}
If $\gamma\equiv1$, ${\mathcal F}_\gamma(X)$ is the usual $n$-dimensional volume of the hypersurface $X$ 
(that is the $n$-dimensional volume of $M$ with the metric induced by $X$).

\subsection{Wulff shape and convexity of integrands} \label{con}

\begin{definition}\label{supp}{\rm 
Assume that $S$ is a closed hypersurface in ${\mathbb R}^{n+1}$ that is the boundary of a bounded connected open set $\Omega$. 
Denote by $\overline{\Omega}$ the closure of $\Omega$, that is, $\overline{\Omega}=\Omega \cup S$.

(i) $S$ is said to be strictly convex if, for any straight line segment $PQ$ connecting two distinct points $P$ and $Q$ in $S$, $PQ \subset \overline{\Omega}$ and $PQ\cap S=\{P, Q\}$ holds.

(ii) $S$ is said to be convex if, for any straight line segment $PQ$ connecting two points $P$ and $Q$ in $S$, $PQ \subset \overline{\Omega}$ holds.
}\end{definition}

A positive continuous function $\gamma:S^n \to {\mathbb R}_{>0}$ is called a convex integrand (or, simply, convex) if 
the set $\{\gamma(\nu)^{-1}\nu \;|\; \nu \in S^n\}$ is convex. This condition is equivalent to the condition that the homogeneous extension $\overline{\gamma}:{\mathbb R}^{n+1} \to {\mathbb R}_{\ge 0}$ of $\gamma$ that is defined as 
\begin{equation}\label{ext}
\overline{\gamma}(rX) := r\gamma(X), \ \forall X \in S^n, \ \forall r \ge 0,
\end{equation}
is a convex function (that is, $\overline{\gamma}(X+Y)\le \overline{\gamma}(X)+\overline{\gamma}(Y)$, $\forall X, Y \in {\mathbb R}^{n+1}$) (cf. \cite{T}).

Now we assume that $\gamma:S^n \rightarrow {\mathbb R}_{>0}$ is a positive continuous function. 
Such a function is sometimes called an integrand. 
The boundary $W_\gamma$ of the convex set 
$
\tilde{W}[\gamma] :=\cap _{\nu \in S^n} 
\bigl\{
X \in {\mathbb R}^{n+1} |\langle X, \nu\rangle \le \gamma(\nu)
\bigr\}
$ 
is called the Wulff shape for $\gamma$, where $\langle \ , \ \rangle$ stands for the standard inner product in ${\mathbb R}^{n+1}$. 
We should remark that originally $\tilde{W}[\gamma]$ itself was called the Wulff shape. However, in this paper, because we mainly study hypersurfaces, we call $W_\gamma:=\partial \tilde{W}[\gamma]$ the Wulff shape. 
In the special case where $\gamma\equiv 1$, $W_\gamma$ coincides with $S^n$.

$W_\gamma$ is not smooth in general. $W_\gamma$ is smooth and strictly convex 
if and only if $\gamma$ is of $C^2$ and the $n\times n$ matrix $D^2\gamma+\gamma\cdot 1$ is positive definite at any point in $S^n$, where $D^2\gamma$ is the Hessian of $\gamma$ on $S^n$ and $1$ is the identity matrix of size $n$. 

\subsection{The Cahn-Hoffman map and anisotropic mean curvature} \label{CHM}

From now on, we assume that $\gamma:S^n \to {\mathbb R}_{>0}$ is a positive $C^2$ function. 
We use the same notations as in \S \ref{pw}. 

The Cahn-Hoffman map  $\xi_\gamma:S^n \to {\mathbb R}^{n+1}$ (for $\gamma$) is defined as
\begin{equation}\label{CHE}
\xi(\nu) :=\xi_\gamma(\nu):=D\gamma+\gamma(\nu)\nu, \quad \nu \in S^n,
\end{equation} 
here $T_\nu(S^n)$ for each $\nu \in S^n$ is naturally identified with a hyperplane in ${\mathbb R}^{n+1}$. 
The Wulff shape $W_\gamma$ is a subset of the image $\hat{W}_\gamma:=\xi_\gamma(S^n)$ of $\xi_\gamma$.  
$\hat{W}_\gamma=W_\gamma$ holds if and only if $\gamma$ is convex, and in this case $\gamma$ is the support function of $W_\gamma$ (cf. \cite{K2018}).

A simple calculation shows that $\xi_\gamma$ is represented by using the homogeneous extension $\overline{\gamma}:{\mathbb R}^{n+1} \to {\mathbb R}_{\ge 0}$ of $\gamma$ as
\begin{equation}\label{homoge2}
\xi_\gamma(\nu)=\overline D\overline\gamma|_{\nu}, \quad \nu \in S^n,
\end{equation}
where $\overline D$ is the gradient in ${\mathbb R}^{n+1}$.

Let $X:M=\cup _{i=1}^k M_i \rightarrow{\mathbb R}^{n+1}$ be a piecewise-$C^2$ weak immersion. 
The Cahn-Hoffman field $\tilde{\xi}_i$ along $X_i=X|_{M_i}$ for $\gamma$ 
(or the anisotropic Gauss map of $X$ for $\gamma$
) is defined as
$
\tilde{\xi}_i:=\xi_\gamma\circ \nu_i:M_i \to {\mathbb R}^{n+1}$. 
The linear map $S^\gamma_p:T_p M_i \to T_p M_i$ given by the $n\times n$ matrix 
$
S^\gamma:=-d\tilde{\xi}_i
$
is called the anisotropic shape operatior of $X_i$.

\begin{definition}[anisotropic curvatures, cf. \cite{R}, \cite{HLMG2009}]\label{curv}{\rm 
(i) The eigenvalues of $S^\gamma$ are called the anisotropic principal
curvatures of $X$. We denote them by $k^\gamma_1, \cdots, k^\gamma_n$.

(ii) Let $\sigma^\gamma_r$ be the elementary symmetric functions of $k^\gamma_1, \cdots, k^\gamma_n$, that is
\begin{equation}\label{curv2}
\sigma^\gamma_r :=\sum_{1\le l_1< \cdots <l_r\le n} k^\gamma_{l_1}\cdots k^\gamma_{l_r}, \quad r=1, \cdots, n.
\end{equation}
We set $\sigma^\gamma_0:=1$. 
The $r$-th anisotropic mean curvature $H^\gamma_r$ of $X$ is defined by $H^\gamma_r:=\sigma^\gamma_r/{}_nC_r$.

(iii) $\Lambda^\gamma:=H^\gamma_1$ is called  the anisotropic mean curvature of $X$. $\Lambda^\gamma$ is often written as $\Lambda$ for simplicity.  
}\end{definition}

Recall that $S(X)$ is the set of all singular points of $X$. For any $p\in M_i\setminus S(X)$, we may write $\tilde{\xi}(p):=\tilde{\xi}_\gamma(p):=\tilde{\xi}_i(p)$. 


\begin{fact}[Euler-Lagrange equations, Koiso \cite{K2018}. For $n=2$, see B. Palmer \cite{P2012}]\label{EL} 
A piecewise-$C^2$ weak immersion $\displaystyle X:M=\cup_{i=1}^k M_i \rightarrow{\mathbb R}^{n+1}$ 
is a critical point of the anisotropic energy $\displaystyle {\mathcal F}_\gamma(X)=\int_M \gamma(\nu) \:dA$ for volume-preserving variations if and only if 

(i) \ The anisotropic mean curvature $\Lambda$ of $X$ is constant on $M\setminus S(X)$, and 

(ii) $\tilde{\xi}_i(\zeta)-\tilde{\xi}_j(\zeta) \in T_\zeta M_i\cap T_\zeta M_j =T_\zeta(\partial M_i\cap \partial M_j)$ holds at any  $\zeta \in \partial M_i\cap \partial M_j$, where a tangent space of a submanifold of ${\mathbb R}^{n+1}$ is naturally identified with a linear subspace of ${\mathbb R}^{n+1}$. 
\end{fact}

\begin{definition}[\cite{K2018}]\label{defR}{\rm
A piecewise-$C^2$ weak immersion $\displaystyle X:M=\cup_{i=1}^k M_i \rightarrow{\mathbb R}^{n+1}$ is called a hypersurface with constant anisotropic mean curvature (CAMC) if both of (i) and (ii) in Fact \ref{EL} hold. 
}\end{definition}

\begin{fact}[\cite{KP2005}, \cite{K2018}]\label{CHR}{\rm 
Since $\xi^{-1}$ gives the unit normal vector field $\nu_{\xi}$ for the Cahn-Hoffman map $\xi$, and since $
\Lambda=-(1/n){\rm trace}_M\;d\tilde\xi
$ holds, we have the following.

(i) The anisotropic mean curvature of the Cahn-Hoffman map $\xi:S^n \to {\mathbb R}^{n+1}$ is $-1$ at any regular point $\nu \in S^n$ with respect to the unit normal $\nu$.  
Hence, particularly the anisotropic mean curvature of the Wulff shape (for the outward-pointing unit normal) is $-1$ at any regular point. 

(ii) Set $A:=D^2\gamma+\gamma\cdot 1$. Take the unit normal vector $\tilde{\nu}$ of the 
Cahn-Hoffman map $\xi$ so that $\tilde{\nu}=\nu$ holds at any point $\nu \in S^n$ where $\det A>0$ holds. 
Let $\{e_1, \cdots, e_n\}$ be an orthonormal basis of $T_\nu S^n$ such that $\{\tilde{\nu}, f_1, \cdots, f_n\}$ is compatible with the canonical orientation of ${\mathbb R}^{n+1}$, where $f_j:=(d_\nu\xi)(e_j)$. 
Then, at each point $\nu \in S^n$ satisfying $\det A<0$, $\tilde{\nu}=-\nu$ holds. 
If $\Lambda(\nu)$ is the anisotropic mean curvature of $\xi$ at $\nu$ with respect to $\tilde{\nu}$, then
\begin{equation}
\Lambda(\nu)= \left \{
\begin{array}{l}
-1, \ (\det A(\nu)>0),\\
1, \ (\det A(\nu)<0).
\end{array}
\right.
\end{equation} 
}\end{fact}
%
%
%
%


\section{An example of dimension $1$}\label{ex1}

We give an example which will be used to prove Theorems \ref{nonuni1}, \ref{nonuni2}, and \ref{AMCF} in the following sections. 
Throughout this section $\gamma:S^1\to{\mathbb R}_{>0}$ is  the function defined by 
\begin{equation}\label{exk2017-1}
\gamma(e^{i\theta}):=\cos^6\theta+\sin^6\theta,
\end{equation}
 where ${\mathbb R}^2$ is identified with ${\mathbb C}$. 

\begin{lemma}\label{lem1}
The Cahn-Hoffman map $\xi_\gamma:S^1 \to {\mathbb R}^2$ for $\gamma$ is represented as follows.
\begin{eqnarray}
\xi_\gamma(e^{i\theta})
=&&\!\!\!\!\!\!\!\!
\bigl((\cos\theta)(\cos^6\theta+6\cos^4\theta\sin^2\theta-5\sin^6\theta), \nonumber\\
&&(\sin\theta)(-5\cos^6\theta+6\cos^4\theta\sin^2\theta+\sin^6\theta)
\bigr).\label{cf1}
\end{eqnarray}
In other words, we have
\begin{equation}\label{cf2}
\xi_\gamma(\nu)
=
\bigl(
\nu_1(\nu_1^6+6\nu_1^4\nu_2^2-5\nu_2^6), 
\nu_2(-5\nu_1^6+6\nu_1^2\nu_2^4+\nu_2^6)
\bigr), \quad \nu=(\nu_1, \nu_2) \in S^1.
\end{equation}
Set $\nu_2=\pm \sqrt{1-\nu_1^2}$. Then 
\begin{eqnarray}
\xi_\gamma(\nu)&=&
\Biggl(
-\nu_1(9\nu_1^4-15\nu_1^2+5), 
\mp\sqrt{1-\nu_1^2}(9\nu_1^4-3\nu_1^2-1)
\Biggr) \label{cf4}
\\
&=&
\Biggl(
-(9\cos^4\theta-15\cos^2\theta+5)\cos\theta, 
-(9\cos^4\theta-3\cos^2\theta-1)\sin\theta
\Biggr)
\label{cf7}
\end{eqnarray}
holds. 
\end{lemma}

\noindent{\it Proof.} \ 
We use the formula (\ref{homoge2}), that is
$$
\xi_\gamma(\nu)=\overline D\overline\gamma|_{\nu}, \quad \nu \in S^n.
$$
$\overline\gamma$ is given by
\begin{equation}\label{gmb0}
\overline\gamma(\nu_1, \nu_2)=\frac{\nu_1^6+\nu_2^6}
{(\nu_1^2+\nu_2^2)^{5/2}}.
\end{equation}
Hence, we have
\begin{eqnarray}
\overline\gamma_{\nu_1}
&=&
\frac{\nu_1(\nu_1^6+6\nu_1^4\nu_2^2-5\nu_2^6)}
{(\nu_1^2+\nu_2^2)^{7/2}}, \label{cfx}\\
\overline\gamma_{\nu_2}
&=&
\frac{\nu_2(-5\nu_1^6+6\nu_1^2\nu_2^4+\nu_2^6)}
{(\nu_1^2+\nu_2^2)^{7/2}}.\label{cfy}
\end{eqnarray}
(\ref{homoge2}) with (\ref{cfx}) and (\ref{cfy}) gives
\begin{equation}\label{cf3}
\xi_\gamma(\nu)
=
\bigl(
\nu_1(\nu_1^6+6\nu_1^4\nu_2^2-5\nu_2^6), 
\nu_2(-5\nu_1^6+6\nu_1^2\nu_2^4+\nu_2^6)
\bigr), \quad \nu=(\nu_1, \nu_2) \in S^1,
\end{equation}
which gives the desired formulas (\ref{cf1}), (\ref{cf2}). 
Inserting $\nu_2=\pm \sqrt{1-\nu_1^2}$ to (\ref{cf2}), we obtain (\ref{cf4}), (\ref{cf7}). 
\hfill $\Box$

\vskip0.5truecm
The image $\xi_\gamma(S^1)$ of the Cahn-Hoffman map $\xi_\gamma$ and the Wulff shape $W_\gamma$ are shown in Figure \ref{fig:ex1-1}. 
Lemma \ref{table} below gives rigorous explanations of these shapes. 
\vskip0.5truecm
    \begin{figure}[H]
  \centering   \includegraphics[width=20mm,height=20mm,angle=0]{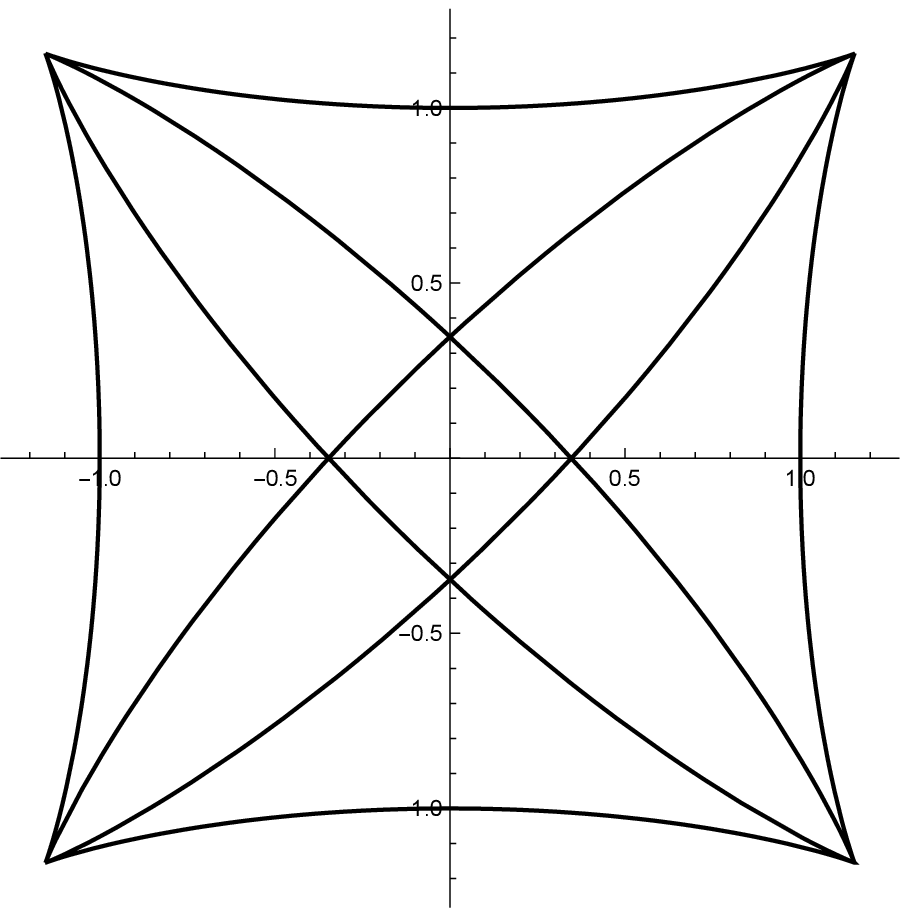}
\hspace{15mm}  \includegraphics[width=10mm,height=10mm,angle=0]{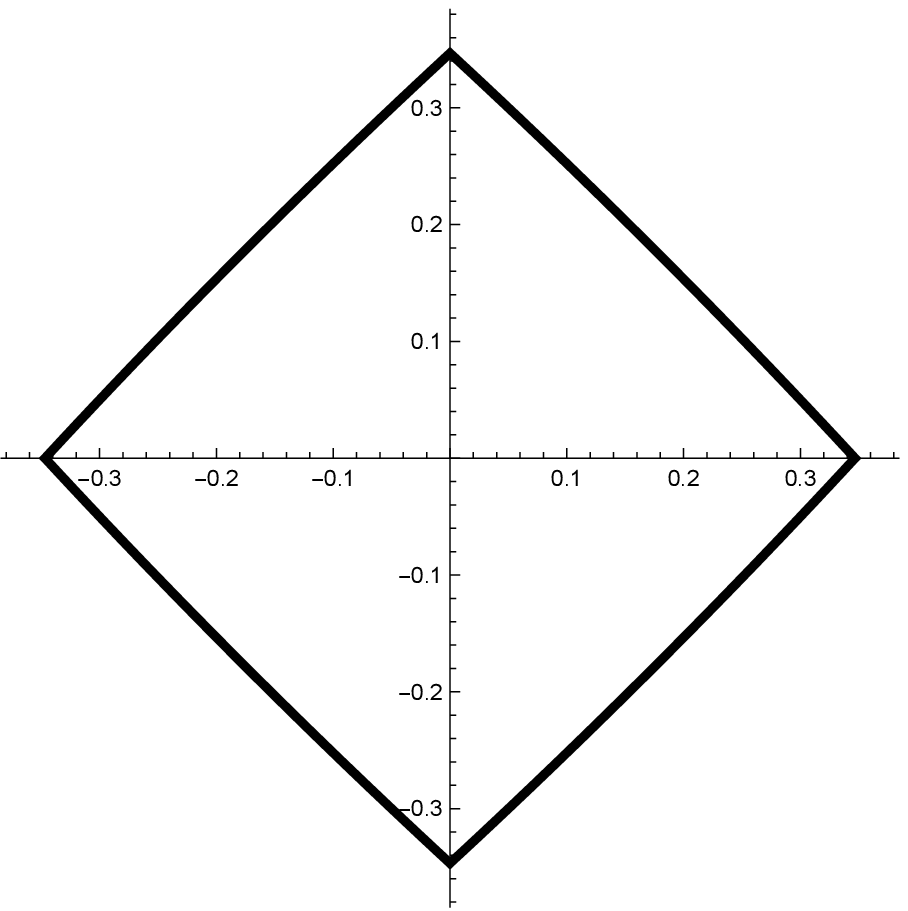} 
 
\caption{The image $\xi_\gamma(S^1)$ of the Cahn-Hoffman map (left) and the Wulff shape $W_\gamma$ (right) for $\gamma$ defined by (\ref{exk2017-1}). $W_\gamma$ is a subset of $\xi_\gamma(S^1)$.}
      \label{fig:ex1-1}
    \end{figure}
    
Now we study the singular points of $\xi_\gamma$. 
\begin{lemma}\label{lem2}
(i) The the set $S(\xi_\gamma)$ of all singular points of $\xi_\gamma$ is given by 
\begin{eqnarray}
S(\xi_\gamma)
&=&\Biggl\{
p_1:=\Bigl(\frac{\sqrt{5}+1}{2\sqrt{3}}, \frac{\sqrt{5}-1}{2\sqrt{3}}\Bigr), \ 
p_2:=\Bigl(\frac{\sqrt{5}-1}{2\sqrt{3}}, \frac{\sqrt{5}+1}{2\sqrt{3}}\Bigr), \nonumber\\
&& p_3:=\Bigl(-\frac{\sqrt{5}-1}{2\sqrt{3}}, \frac{\sqrt{5}+1}{2\sqrt{3}}\Bigr), \ p_4:=\Bigl(-\frac{\sqrt{5}+1}{2\sqrt{3}}, \frac{\sqrt{5}-1}{2\sqrt{3}}\Bigr), 
\ \nonumber\\
&&p_5:=\Bigl(-\frac{\sqrt{5}+1}{2\sqrt{3}}, -\frac{\sqrt{5}-1}{2\sqrt{3}}\Bigr), 
\ p_6:=\Bigl(-\frac{\sqrt{5}-1}{2\sqrt{3}}, -\frac{\sqrt{5}+1}{2\sqrt{3}}\Bigr), \nonumber\\
&&p_7:=\Bigl(\frac{\sqrt{5}-1}{2\sqrt{3}}, -\frac{\sqrt{5}+1}{2\sqrt{3}}\Bigr), \ p_8:=\Bigl(\frac{\sqrt{5}+1}{2\sqrt{3}}, -\frac{\sqrt{5}-1}{2\sqrt{3}}\Bigr)
\Biggr\}. \label{sin1}
\end{eqnarray}

(ii) Set 
$$
p_j=(\cos\theta_j, \sin\theta_j), \quad j=1, \cdots, 8,
$$
where $\theta_1$ is chosen so that $0< \theta_1 < \pi/2$ holds. 
Then, by choosing suitable principal values for $\theta_1, \cdots, \theta_8$, one can write
\begin{eqnarray}
\theta_2&=&\pi/2-\theta_1, \ 
\theta_3=\pi/2+\theta_1, \ 
\theta_4=\pi-\theta_1, \ \theta_5=\pi+\theta_1, \ \label{t1}\\
\theta_6&=&3\pi/2-\theta_1, \
\theta_7=3\pi/2+\theta_1, \   \theta_8=-\theta_1. \label{t2}
\end{eqnarray}

(iii)
\begin{eqnarray} 
\xi_\gamma(p_8)&=&\xi_\gamma(p_3)=\frac{2}{\sqrt{3}}(1, 1), \  \xi_\gamma(p_2)=\xi_\gamma(p_5)=\frac{2}{\sqrt{3}}(-1, 1), \\
\xi_\gamma(p_4)&=&\xi_\gamma(p_7)=\frac{2}{\sqrt{3}}(-1, -1), \ \xi_\gamma(p_6)=\xi_\gamma(p_1)=\frac{2}{\sqrt{3}}(1, -1)
\end{eqnarray}
holds. 
\end{lemma}

\begin{remark}\label{theta}{\rm
Choose the principal value for $\theta_1$ so that $0<\theta_1<\pi/2$ is satisfied. Then, by computation using Mathematica ver.11.2.0.0, we get 
\begin{equation}\label{comp}
\theta_1\approx 0.1161397636\pi.
\end{equation}
}
\end{remark}

\noindent{\it Proof of Lemma \ref{lem2}.} \ 
Using (\ref{cf4}), we obtain
\begin{equation}\label{dif1}
\xi_\gamma'(\nu_1)
=\Bigl(
-5(9\nu_1^4-9\nu_1^2+1), 
\pm\frac{5\nu_1}{\sqrt{1-\nu_1^2}}
(9\nu_1^4-9\nu_1^2+1)
\Bigr).
\end{equation}
Hence, $(\nu_1, \nu_2)\in S^1$ is a singular point of $\xi_\gamma$ if and only if
\begin{equation}\label{cf5}
9\nu_1^4-9\nu_1^2+1=0
\end{equation}
holds. (\ref{cf5}) is equivalent to 
\begin{equation}\label{cf6}
\nu_1=\frac{\sqrt{5}\pm 1}
{2\sqrt{3}}, \ - \frac{\sqrt{5}\pm 1}
{2\sqrt{3}},
\end{equation}
which proves (i). 

(ii) is a consequence of (i). 

(iii) is obtained by inserting (\ref{sin1}) to (\ref{cf4}). 
\hfill $\Box$

\vskip0.5truecm

Set
$$
\xi_\gamma(\theta)=((\xi_\gamma)_x(\theta), (\xi_\gamma)_y(\theta)):=\xi_\gamma(e^{i\theta}).
$$

\begin{lemma}\label{table}
Increase and decrease of $(\xi_\gamma)_x(\theta)$, $(\xi_\gamma)_y(\theta)$, and the sign of $A=D^2\gamma+\gamma\cdot 1$ are given by the following table. 

\begin{table}[H]
\label{tab1}
\begin{center}
\scalebox{0.86}{
\begin{tabular}{@{}l@{}}
\noalign{\hrule height0.8pt}
\begin{tabular}{|c|c|c|c|c|c|c|c|c|c|c|c|c|c|}
{$\theta$}& {$\theta_8$} & {} & {$0$}& {}& {$\theta_1$} & {} & {$\theta_2$} & {}&{$\pi/2$} & {} & {$\theta_3$} & {}&{$\theta_4$}\\
\hline
{$(\xi_\gamma)'_x(\theta)$}& {$0$} & {$-$} & {$0$}& {$+$}& {$0$} & {$-$} & {$0$} & {$+$}&{$5$} & {$+$} & {$0$} & {$-$}&{$0$}\\
\hline
{$(\xi_\gamma)'_y(\theta)$}& {$0$} & {$-$} & {$-5$}& {$-$}& {$0$} & {$+$} & {$0$} & {$-$}&{$0$} & {$+$} & {$0$} & {$-$}&{$0$}\\
\hline
{$(\xi_\gamma)_x(\theta)$}& {$2/\sqrt{3}$} & {$\searrow$} & {1}& {$\nearrow$}& {$2/\sqrt{3}$} & {$\searrow$} & {$-2/\sqrt{3}$} & {$\nearrow$}
&{$0$} & {$\nearrow$} & {$2/\sqrt{3}$} & {$\searrow$}&{$-2/\sqrt{3}$}\\
\hline
{$(\xi_\gamma)_y(\theta)$}& {$2/\sqrt{3}$} & {$\searrow$} & {$0$} &{$\searrow$}& {$-2/\sqrt{3}$} & {$\nearrow$} & {$2/\sqrt{3}$} & {$\searrow$}
&{$1$} & {$\nearrow$} & {$2/\sqrt{3}$} & {$\searrow$}&{$-2/\sqrt{3}$}\\
\hline
{$A=D^2\gamma+\gamma\cdot 1$}& {$0$} & {$-$} & {$-$} &{$-$}& {$0$} & {$+$} & {$0$} & {$-$}
&{$-$} & {$-$} & {$0$} & {$+$}&{$0$}\\

\hline
\end{tabular}\\
\noalign{\hrule height0.8pt} 
\end{tabular}
}
\end{center}
\end{table}

\begin{table}[H]
\label{tab2}
\begin{center}
\scalebox{0.86}{
\begin{tabular}{@{}l@{}}
\noalign{\hrule height0.8pt}
\begin{tabular}{|c|c|c|c|c|c|c|c|c|c|c|c|c|c|}
{$\theta$}& {$\theta_4$} & {} & {$\pi$}& {}& {$\theta_5$} & {} & {$\theta_6$} & {}&{$3\pi/2$} & {} & {$\theta_7$} & {}&{$\theta_8+2\pi$}\\
\hline
{$(\xi_\gamma)'_x(\theta)$}& {$0$} & {$+$} & {$0$}& {$-$}& {$0$} & {$+$} & {$0$} & {$-$}&{$-5$} & {$-$} & {$0$} & {$+$}&{$0$}\\
\hline
{$(\xi_\gamma)'_y(\theta)$}& {$0$} & {$+$} & {$5$}& {$+$}& {$0$} & {$-$} & {$0$} & {$+$}&{$0$} & {$-$} & {$0$} & {$+$}&{$0$}\\
\hline
{$(\xi_\gamma)_x(\theta)$}& {$-2/\sqrt{3}$} & {$\nearrow$} & {$-1$}& {$\searrow$}& {$-2/\sqrt{3}$} & {$\nearrow$} & {$2/\sqrt{3}$} & {$\searrow$}
&{$0$} & {$\searrow$} & {$-2/\sqrt{3}$} & {$\nearrow$}&{$2/\sqrt{3}$}\\
\hline
{$(\xi_\gamma)_y(\theta)$}& {$-2/\sqrt{3}$} & {$\nearrow$} & {$0$} &{$\nearrow$}& {$2/\sqrt{3}$} & {$\searrow$} & {$-2/\sqrt{3}$} & {$\nearrow$}
&{$-1$} & {$\searrow$} & {$-2/\sqrt{3}$} & {$\nearrow$}&{$2/\sqrt{3}$}\\
\hline
{$A=D^2\gamma+\gamma\cdot 1$}& {$0$} & {$-$} & {$-$} &{$-$}& {$0$} & {$+$} & {$0$} & {$-$}
&{$-$} & {$-$} & {$0$} & {$+$}&{$0$}\\

\hline
\end{tabular}\\
\noalign{\hrule height0.8pt} 
\end{tabular}
}
\end{center}
\end{table}

\end{lemma}

\noindent{\it Proof.} \ 
Using (\ref{cf7}), we obtain
\begin{equation}\label{dif2}
\xi_\gamma'(\theta)
=((\xi_\gamma)'_x(\theta), (\xi_\gamma)'_y(\theta))
=\Bigl(
5(9\cos^4\theta-9\cos^2\theta+1)\sin\theta, 
-5(9\cos^4\theta-9\cos^2\theta+1)\cos\theta
\Bigr).
\end{equation}
For $\gamma(\theta):=\gamma(e^{i\theta})=\cos^6\theta+\sin^6\theta$, we compute
$$
\gamma'(\theta)=-6\cos^5\theta\sin\theta+6\sin^5\theta\cos\theta,
$$
$$
\gamma''(\theta)
=
6(5\cos^2\theta\sin^2\theta-\cos^6\theta-\sin^6\theta),
$$
and hence
\begin{equation}\label{A1}
A=\gamma''(\theta)+\gamma(\theta)
=
-5(9\cos^4\theta-9\cos^2\theta+1)
=-(\xi_\gamma)_x'(\theta)\cdot(\sin\theta)^{-1}
\end{equation}
holds, here the last equality is valid only when $\sin\theta\ne 0$.
These observations with Lemmas \ref{lem1}, \ref{lem2} give the desired result. 
\hfill $\Box$

\vskip0.5truecm
Using Lemma \ref{lem1}, by simple calculations, we obtain the following result. 
\begin{lemma}\label{inter}
Let $\rho_1, \rho_2 \in (0, \pi/2)$ be the solutions of 
\begin{equation}\label{rho12}
\cos\rho_1=\sqrt{
\frac{1+\sqrt{5}}
{6}}, \ 
\cos\rho_2=\sqrt{
\frac{5-\sqrt{5}}
{6}},
\end{equation}
respectively. Then

(i) 
\begin{equation}\label{tr}
0<\theta_1<\rho_1<\pi/4<\rho_2<\theta_2<\pi/2
\end{equation}
holds. 

(ii) Set
\begin{equation}
\alpha:=\sqrt{
\frac{2+2\sqrt{5}}
{3}
}\;
(\sqrt{5}-2).
\end{equation}
The ``inner self-intersection points'' of $\xi_\gamma$ are the following four points. 
\begin{eqnarray}\label{int}
Q_1\!\!\!&:=&\!\!\!(\alpha, 0)=\xi_\gamma(\rho_1)=\xi_\gamma(\rho_2+3\pi/2), \ 
Q_2:=(0, \alpha)=\xi_\gamma(\rho_2)=\xi_\gamma(\rho_1+\pi/2), \nonumber\\
Q_3\!\!\!&:=&\!\!\!(-\alpha, 0)=\xi_\gamma(\rho_1+\pi)=\xi_\gamma(\rho_2+\pi/2), \ 
Q_4:=(0, -\alpha)=\xi_\gamma(\rho_2+\pi)=\xi_\gamma(\rho_1+3\pi/2). \nonumber
\end{eqnarray}
$Q_1, Q_2, Q_3, Q_4$ are the vertices of $W_\gamma$ (Figure \ref{fig:ex1-1}, right). 
\end{lemma}

\begin{remark}\label{theta1}{\rm
By computation using Mathematica ver.11.2.0.0, we get 
\begin{equation}\label{comp2}
\rho_1\approx 0.2374632441\pi, \quad  \rho_2=\pi/2-\rho_1\approx 0.2625367559\pi, \quad 
\alpha\approx 0.3467370642.
\end{equation}
}
\end{remark}

\begin{figure}[H]
 \centering
  \begin{minipage}[b]
  {0.20\linewidth}    
\centering \includegraphics[width=30mm,height=30mm,angle=0]{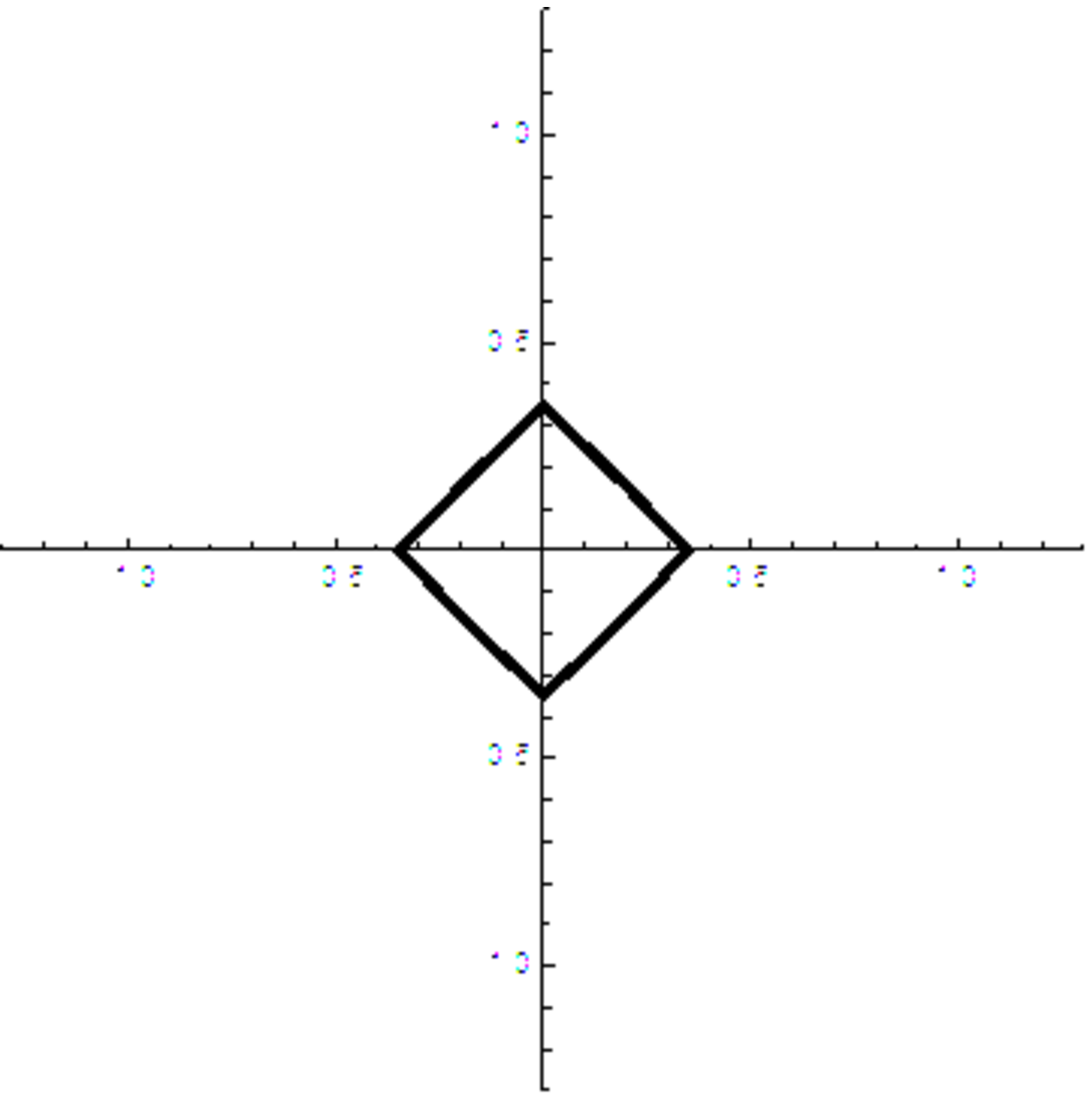}
    \subcaption{Wulff shape $W_\gamma$}\label{fig-a}
  \end{minipage}
  \hspace{5mm}
    \begin{minipage}[b]{0.25\linewidth}
    \centering
    \includegraphics[width=25mm,height=25mm,angle=0]{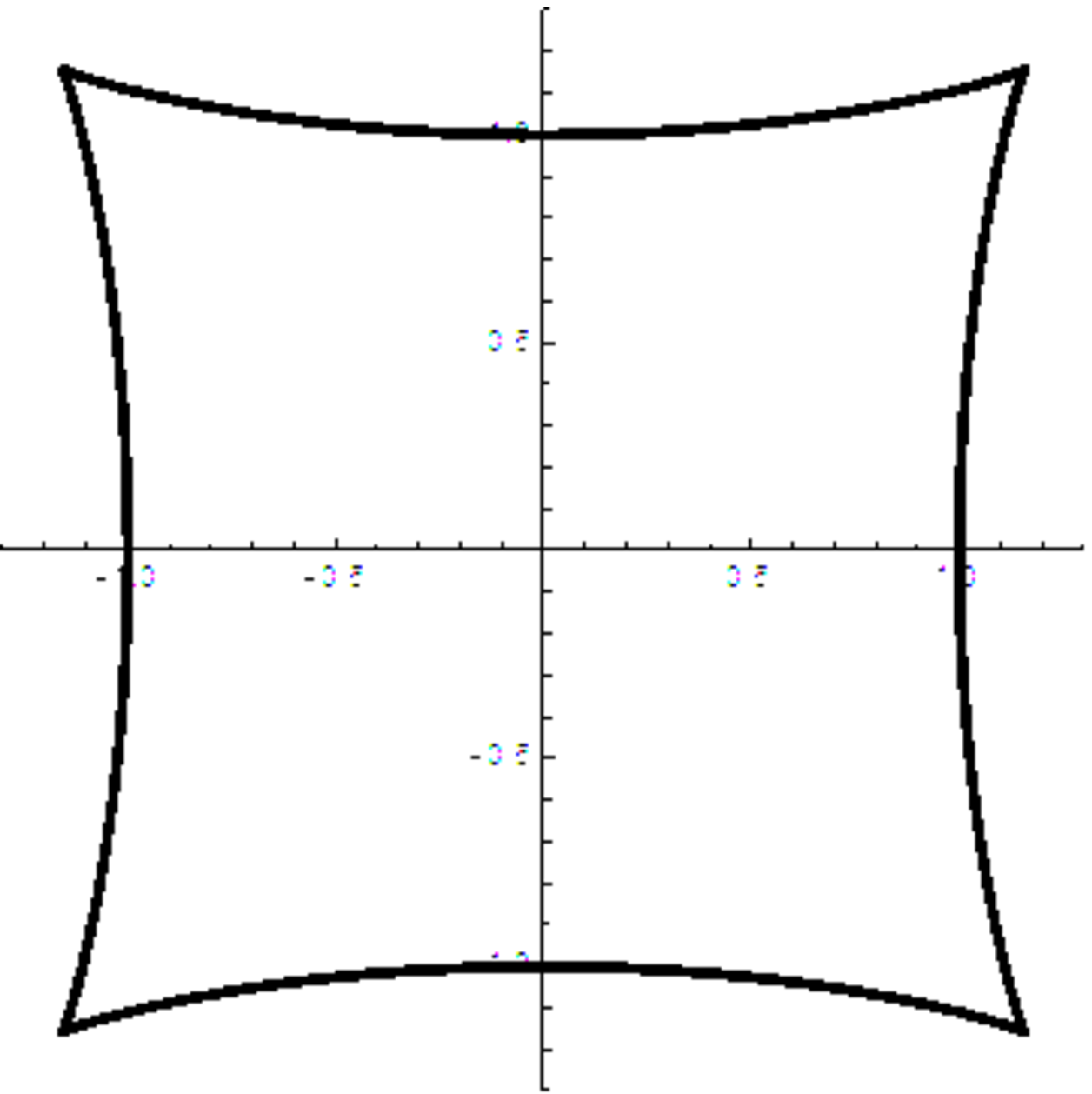}
    \subcaption{CAMC curve $(C_\gamma)_1$}\label{fig-b}
  \end{minipage} 
    \hspace{5mm}
    \begin{minipage}[b]{0.25\linewidth}
    \centering
    \includegraphics[width=25mm,height=25mm,angle=0]{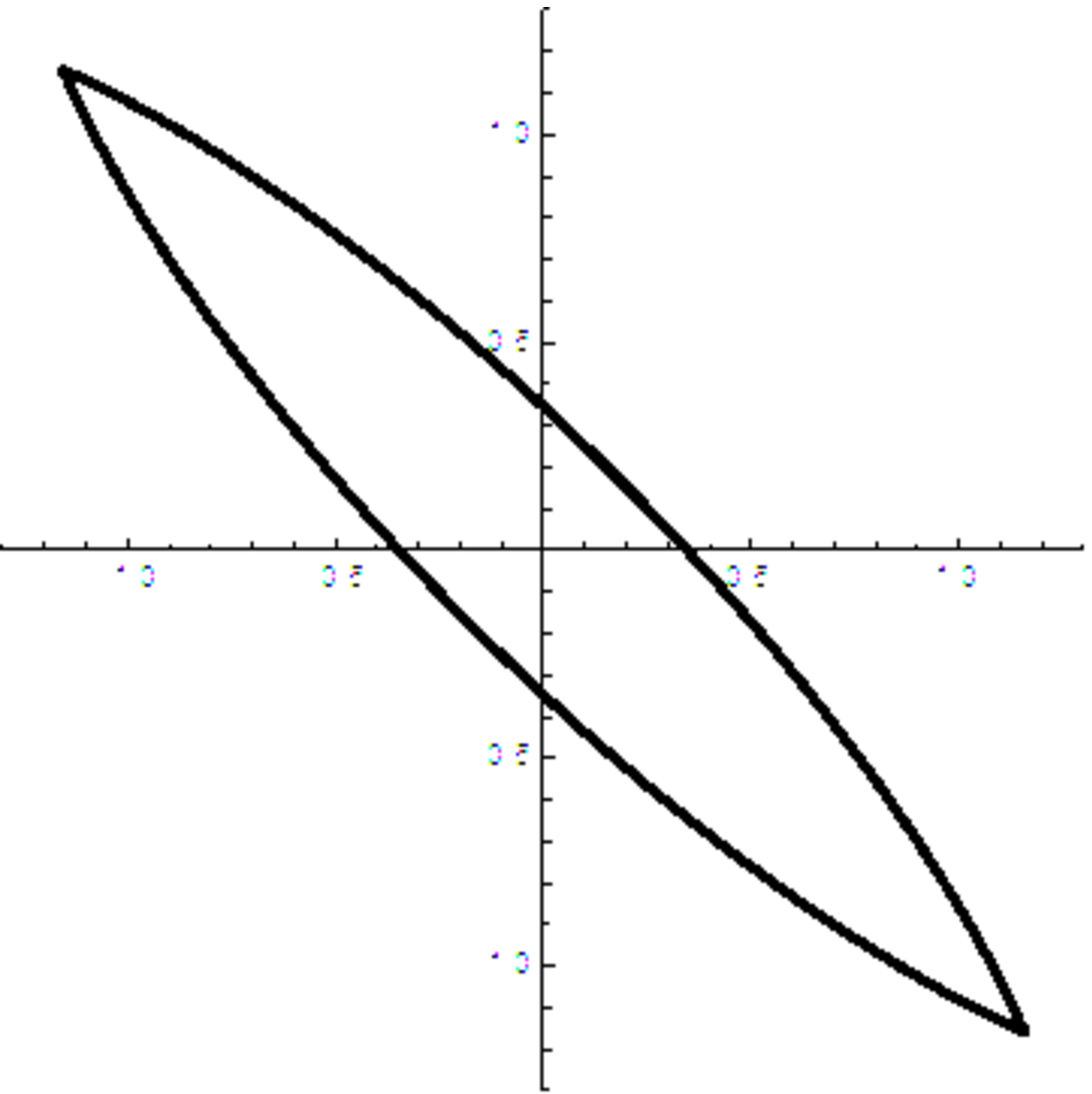}
    \subcaption{CAMC curve $(C_\gamma)_2$}\label{fig-c}
  \end{minipage}
  \hspace{5mm} 
    \begin{minipage}[b]{0.25\linewidth}
    \centering    \includegraphics[width=24mm,height=24mm,angle=0]{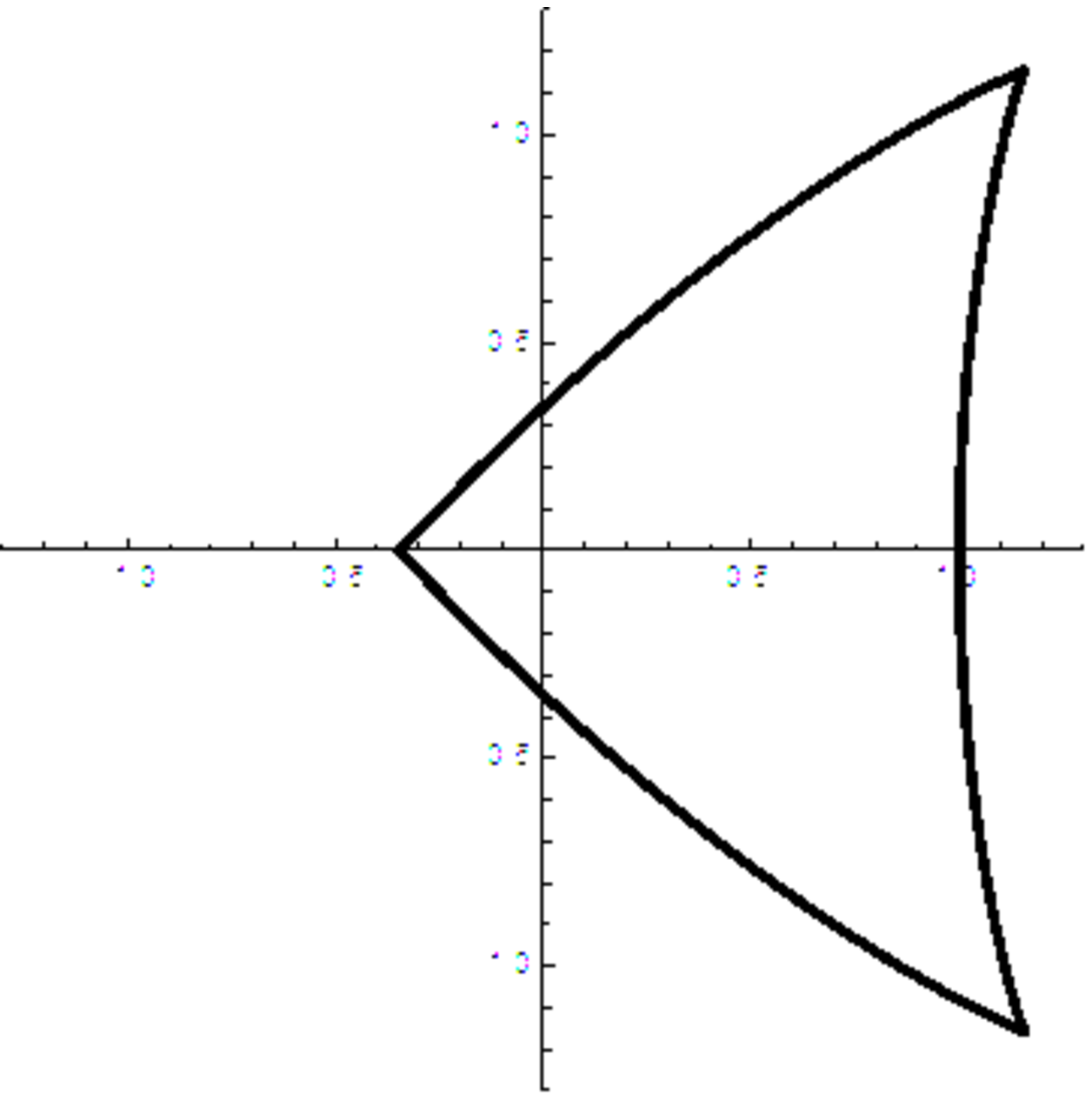}
    \subcaption{CAMC curve $(C_\gamma)_3$}\label{fig-d}
  \end{minipage}   
\hspace{5mm} 
    \begin{minipage}[b]{0.25\linewidth}
    \centering    \includegraphics[width=24mm,height=24mm,angle=0]{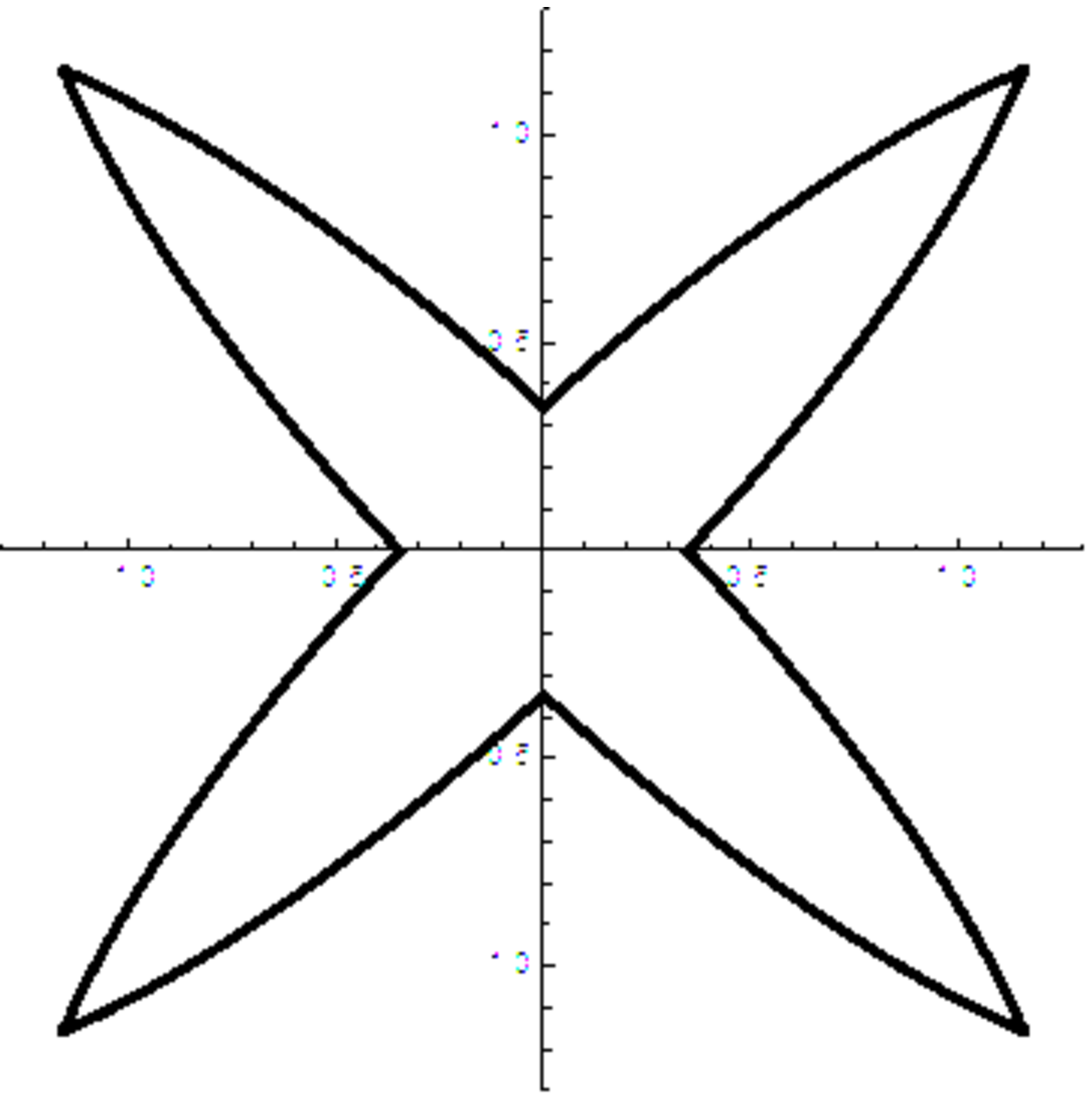}
    \subcaption{CAMC curve $(C_\gamma)_4$}\label{fig-e}
  \end{minipage}   
\hspace{5mm} 
    \begin{minipage}[b]{0.30\linewidth}
    \centering      \includegraphics[width=25mm,height=24mm,angle=0]{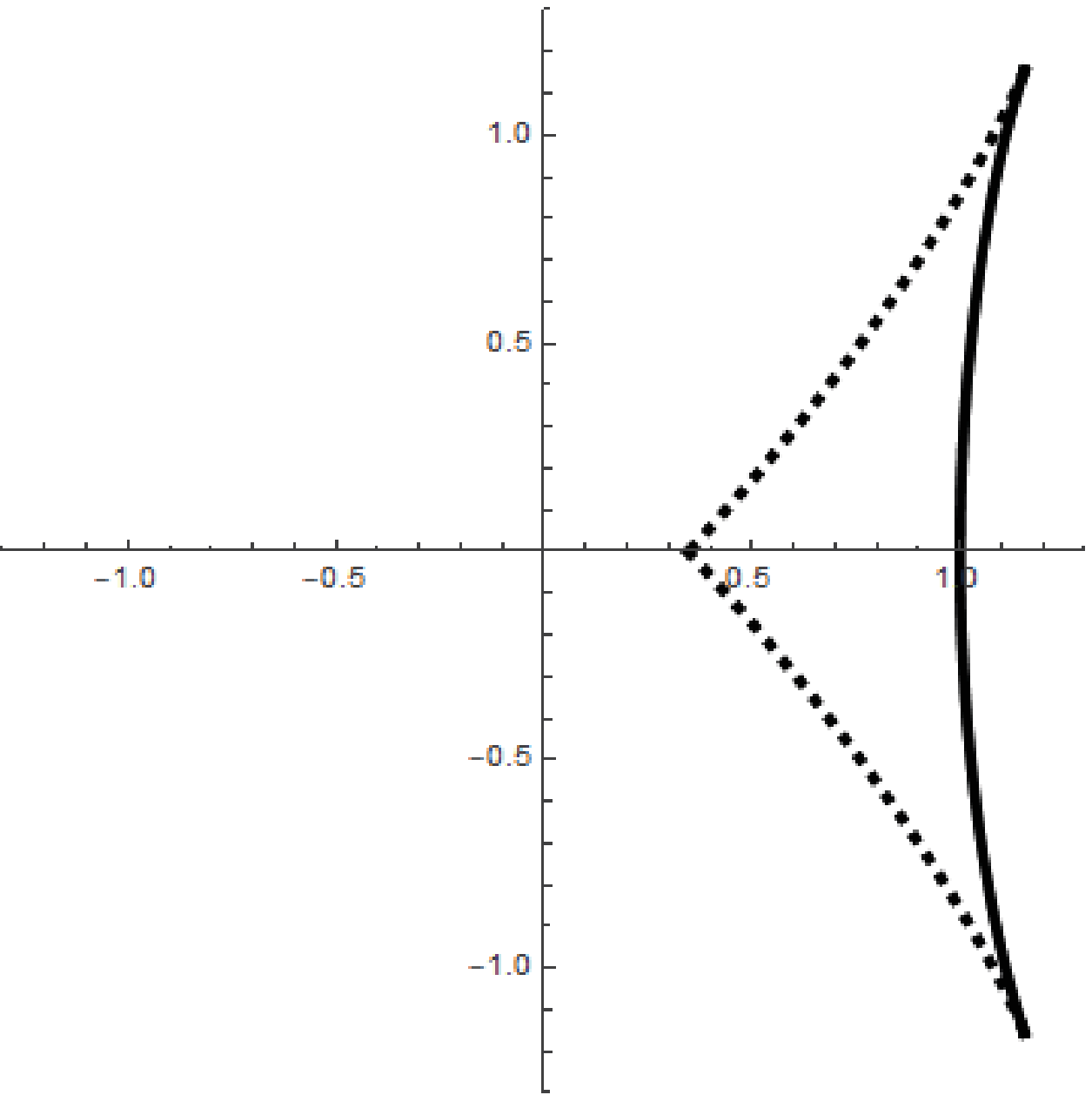} 
        \subcaption{curve $(C_\gamma)_5$. not CAMC}\label{fig-f}
  \end{minipage} 
  \caption{Some of the closed curves  which are subsets of $\xi_\gamma(S^1)$ for $\gamma$ defined by (\ref{exk2017-1}) (cf. Figure \ref{fig:ex1-1}). (a), (b), (c), (d), (e) : The anisotropic (mean) curvature for the outward-pointing normal is $-1$. (f) : For the outward-pointing normal, the anisotropic (mean) curvature is $-1$ at each point in the solid curve, while it is $1$ at each point in the dashed curves. Hence, this curve is not CAMC.}\label{fig-p}
\end{figure}

By using the previous lemmas \ref{lem1}, \ref{lem2}, \ref{table}, and \ref{inter} with Fact \ref{CHR} (i), we obtain the following:
\begin{proposition}\label{prop}
Let   
$\gamma:S^1\to{\mathbb R}_{>0}$ be the function defined by 
\begin{equation}\label{gam0}
\gamma(e^{i\theta}):=\cos^6\theta+\sin^6\theta.
\end{equation}
Let $\theta_1, \theta_2$ be the constants defined in Lemma \ref{lem2}, and let $\rho_1, \rho_2$ be the constants defined in Lemma\ref{inter}.

(i) The Cahn-Hoffman map $\xi_\gamma:S^1 \to {\mathbb R}^2$ is represented as 
\begin{eqnarray}
\xi_\gamma(e^{i\theta})
=&&\!\!\!\!\!\!\!\!
\bigl((\cos\theta)(\cos^6\theta+6\cos^4\theta\sin^2\theta-5\sin^6\theta), \nonumber\\
&&(\sin\theta)(-5\cos^6\theta+6\cos^4\theta\sin^2\theta+\sin^6\theta)
\bigr) \label{cf0}
\end{eqnarray}
(Figure \ref{fig:ex1-1}).

(ii) 
The Wulff shape $W_\gamma$ is given by
\begin{eqnarray}
W_\gamma
=&&\!\!\!\!\! \xi_\gamma
\Bigl(
\{\rho_1\le \theta \le \rho_2\}
\cup
\{\rho_1+\pi/2\le \theta \le \rho_2+\pi/2\} \nonumber\\
&&\cup
\{\rho_1+\pi\le \theta \le \rho_2+\pi\}
\cup
\{\rho_1+3\pi/2\le \theta \le \rho_2+3\pi/2\}
\Bigr) \label{W0}
\end{eqnarray}
(Figure \ref{fig-a}).

(iii) The following four closed curves which are subsets of $\xi_\gamma(S^1)$ are closed piecewise-$C^\infty$ CAMC curves (Figures \ref{fig-b}, \ref{fig-c}, \ref{fig-d}, \ref{fig-e}). 
\begin{eqnarray}
(C_\gamma)_1
:=&&\!\!\!\!\! \xi_\gamma
\Bigl(
\{\theta_2\le \theta \le \theta_3\}
\cup
\{\theta_4\le \theta \le \theta_5\} \nonumber\\
&&\cup
\{\theta_6\le \theta \le \theta_7\}
\cup
\{\theta_8\le \theta \le \theta_1\}
\Bigr), \label{CA01}
\end{eqnarray}
\begin{equation}\label{CA02}
(C_\gamma)_2
:= \xi_\gamma
\Bigl(
\{\theta_1\le \theta \le \theta_2\}
\cup
\{\theta_5\le \theta \le \theta_6\} 
\Bigr), 
\end{equation}
\begin{equation}\label{CA03}
(C_\gamma)_3
:= \xi_\gamma
\Bigl(
\{\theta_8\le \theta \le \theta_1\}
\cup
\{\theta_3\le \theta \le \rho_2+\pi/2\}
\cup
\{\rho_1+\pi\le \theta \le \theta_6\} 
\Bigr), 
\end{equation}
\begin{eqnarray}
(C_\gamma)_4
:=&&\!\!\!\!\! \xi_\gamma
\Bigl(
\{\theta_1\le \theta \le \rho_1\}
\cup
\{\rho_2\le \theta \le \theta_2\} \nonumber\\
&&\cup
\{\theta_3\le \theta \le \rho_1+\pi/2\}
\cup
\{\rho_2+\pi/2\le \theta \le \theta_4\}
\Bigr)\nonumber\\
&&\cup
\{\theta_5\le \theta \le \rho_1+\pi\}
\cup
\{\rho_2+\pi\le \theta \le \theta_6\}
\Bigr)\nonumber\\
&&\cup
\{\theta_7\le \theta \le \rho_1+3\pi/2\}
\cup
\{\rho_2+3\pi/2 \le \theta \le \theta_8+2\pi\}
\Bigr). \label{CA04}
\end{eqnarray}
The anisotropic (mean) curvature for the outward-pointing normal is $-1$.

(iv) The following closed curve which is a subset of $\xi_\gamma(S^1)$ is a closed piecewise-$C^\infty$ curve (Figure \ref{fig-f}). 
\begin{equation}\label{CA05}
(C_\gamma)_5
:= \xi_\gamma
\Bigl(
\{-\rho_1\le \theta \le \rho_1\}
\Bigr). 
\end{equation}
Its anisotropic (mean) curvature is not constant, because, for the outward-pointing normal, it is $-1$ at each point in the solid curve:
\begin{equation}\label{CA04a}
\xi_\gamma
\Bigl(
\{\theta_8\le \theta \le \theta_1\}
\Bigr), 
\end{equation}
while it is $1$ at each point in the dashed curves:
\begin{equation}\label{CA04B}
\xi_\gamma
\Bigl(
\{-\rho_1\le \theta \le \theta_8\}
\cup 
\{\theta_1\le \theta \le \rho_1\}
\Bigr). 
\end{equation}
\end{proposition}

\begin{remark}\label{rem3}{\rm
Proposition \ref{prop} proves Theorem \ref{nonuni1} for $n=1$. 
}\end{remark}

\begin{figure}[H]
 \centering
  \begin{minipage}[b]
  {0.20\linewidth}    
\centering \includegraphics[width=25mm,height=25mm,angle=0]{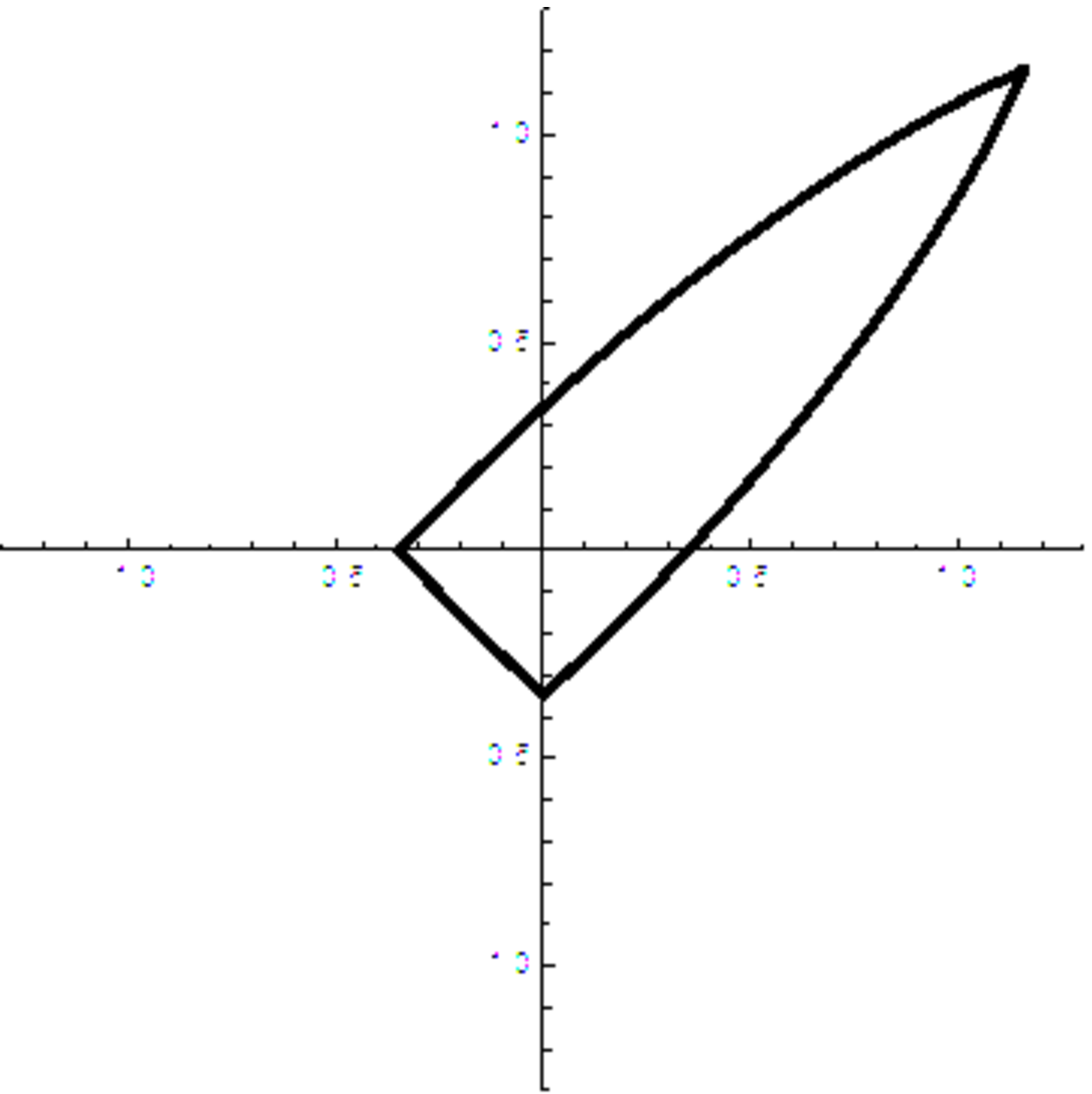}
    \subcaption{CAMC curve}\label{fig-2a}
  \end{minipage}  
  \hspace{5mm}    
    \begin{minipage}[b]{0.25\linewidth}
    \centering   \includegraphics[width=25mm,height=25mm,angle=0]{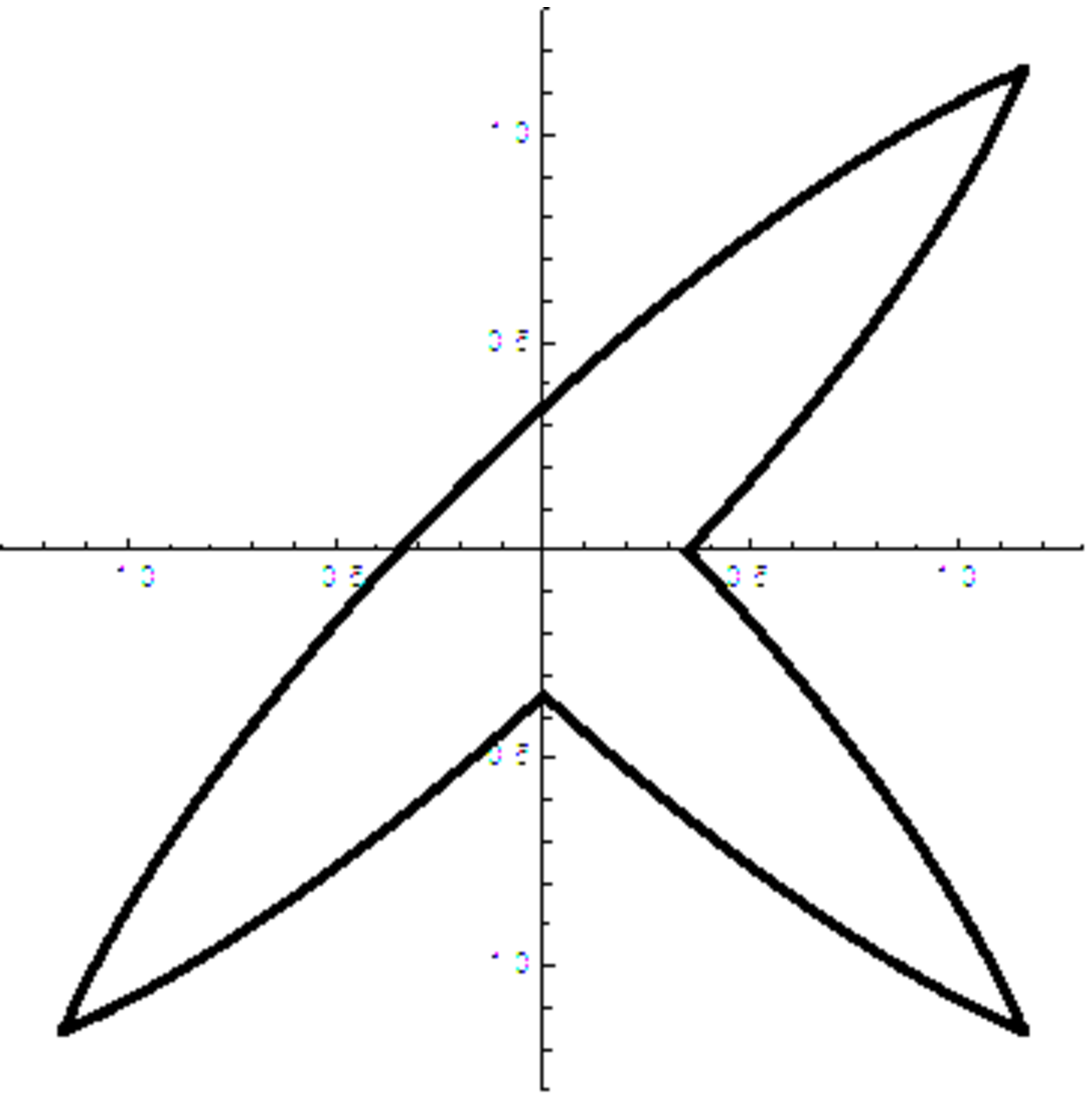}
    \subcaption{CAMC curve}\label{fig-2b}
  \end{minipage}     
    \hspace{5mm}   
    \begin{minipage}[b]{0.25\linewidth}
    \centering    \includegraphics[width=25mm,height=25mm,angle=0]{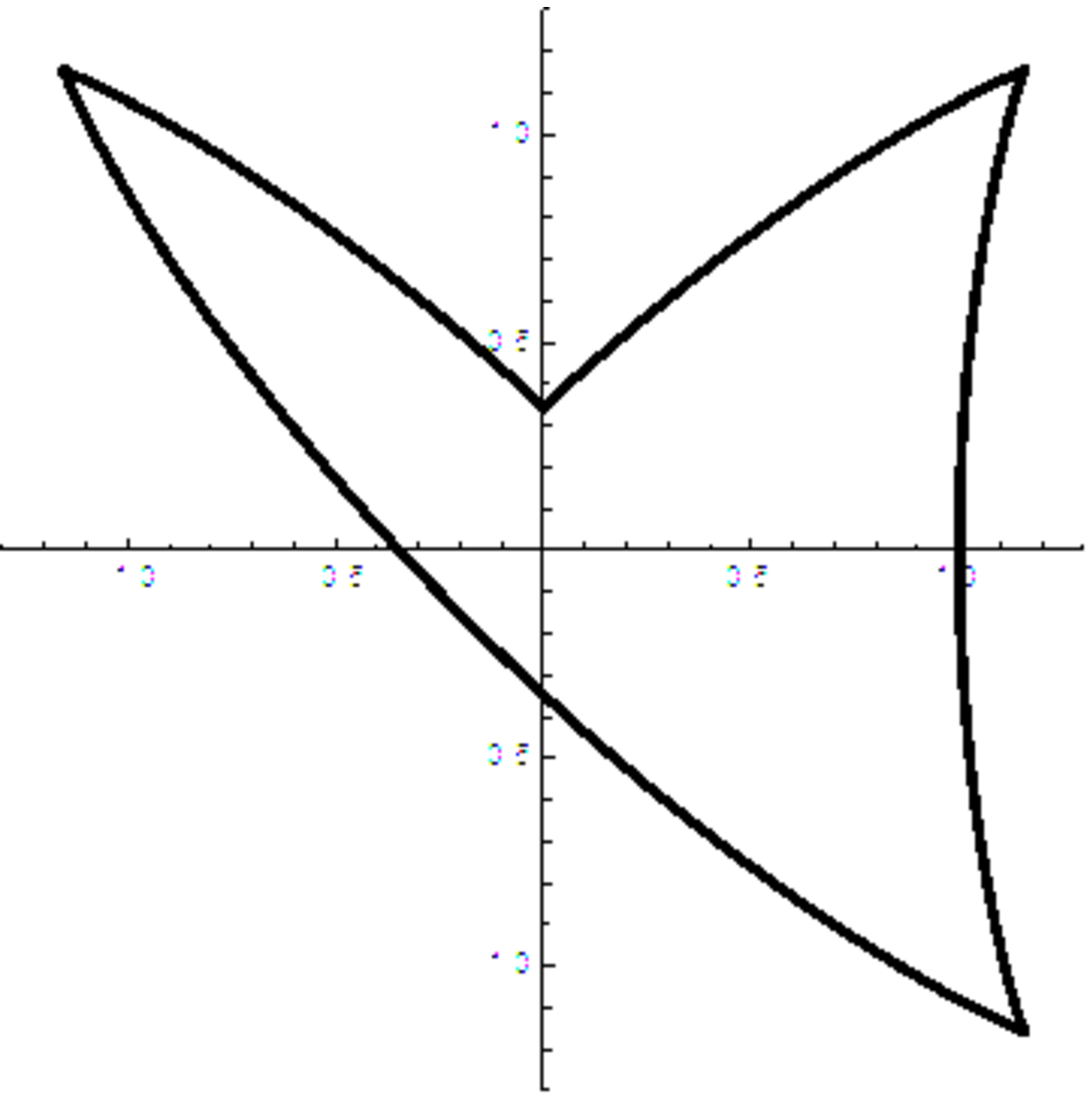}
    \subcaption{CAMC curve}\label{fig-2c}
  \end{minipage}  
  \hspace{5mm}     
    \begin{minipage}[b]{0.25\linewidth}
    \centering    \includegraphics[width=24mm,height=24mm,angle=0]{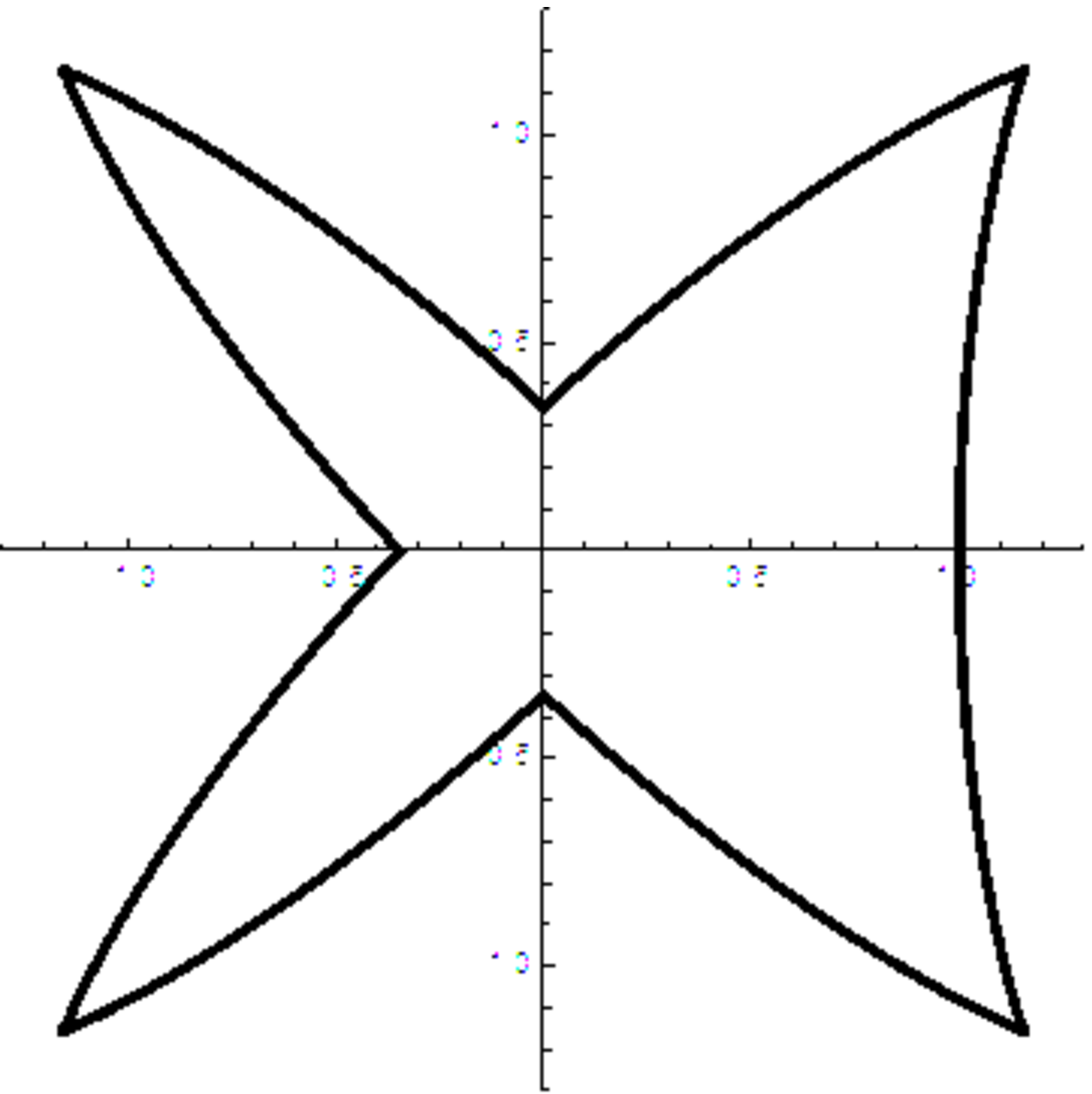}
    \subcaption{CAMC curve}\label{fig-2d}
  \end{minipage}   
\hspace{5mm}     
    \begin{minipage}[b]{0.25\linewidth}
    \centering    \includegraphics[width=24mm,height=24mm,angle=0]{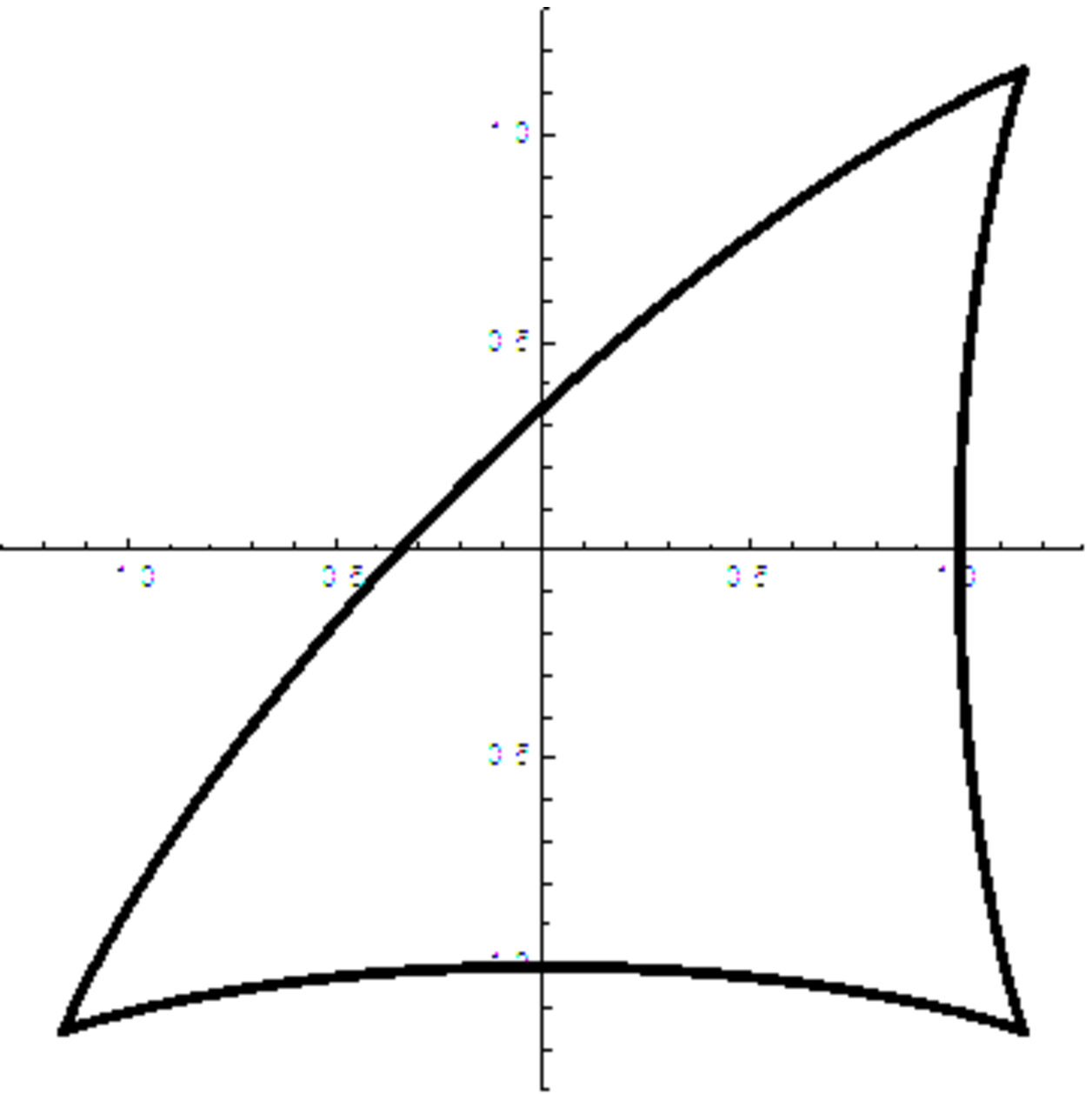}
    \subcaption{CAMC curve}\label{fig-2e}
  \end{minipage}   
\hspace{5mm}     
    \begin{minipage}[b]{0.30\linewidth}
    \centering      \includegraphics[width=25mm,height=24mm,angle=0]{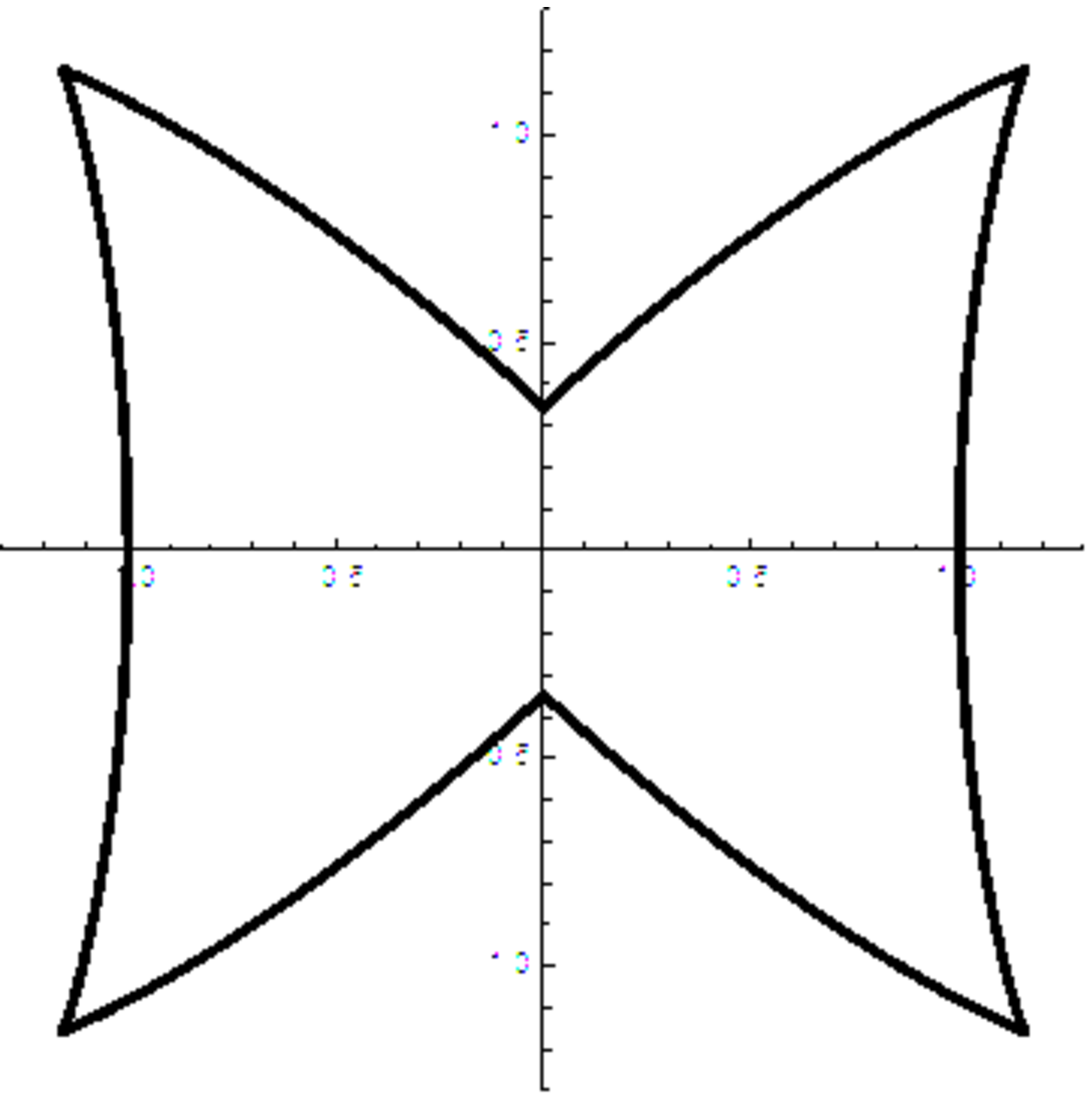} 
        \subcaption{CAMC curve}\label{fig-2f}
  \end{minipage}   
  \hspace{5mm}     
    \begin{minipage}[b]{0.30\linewidth}
    \centering      \includegraphics[width=25mm,height=24mm,angle=0]{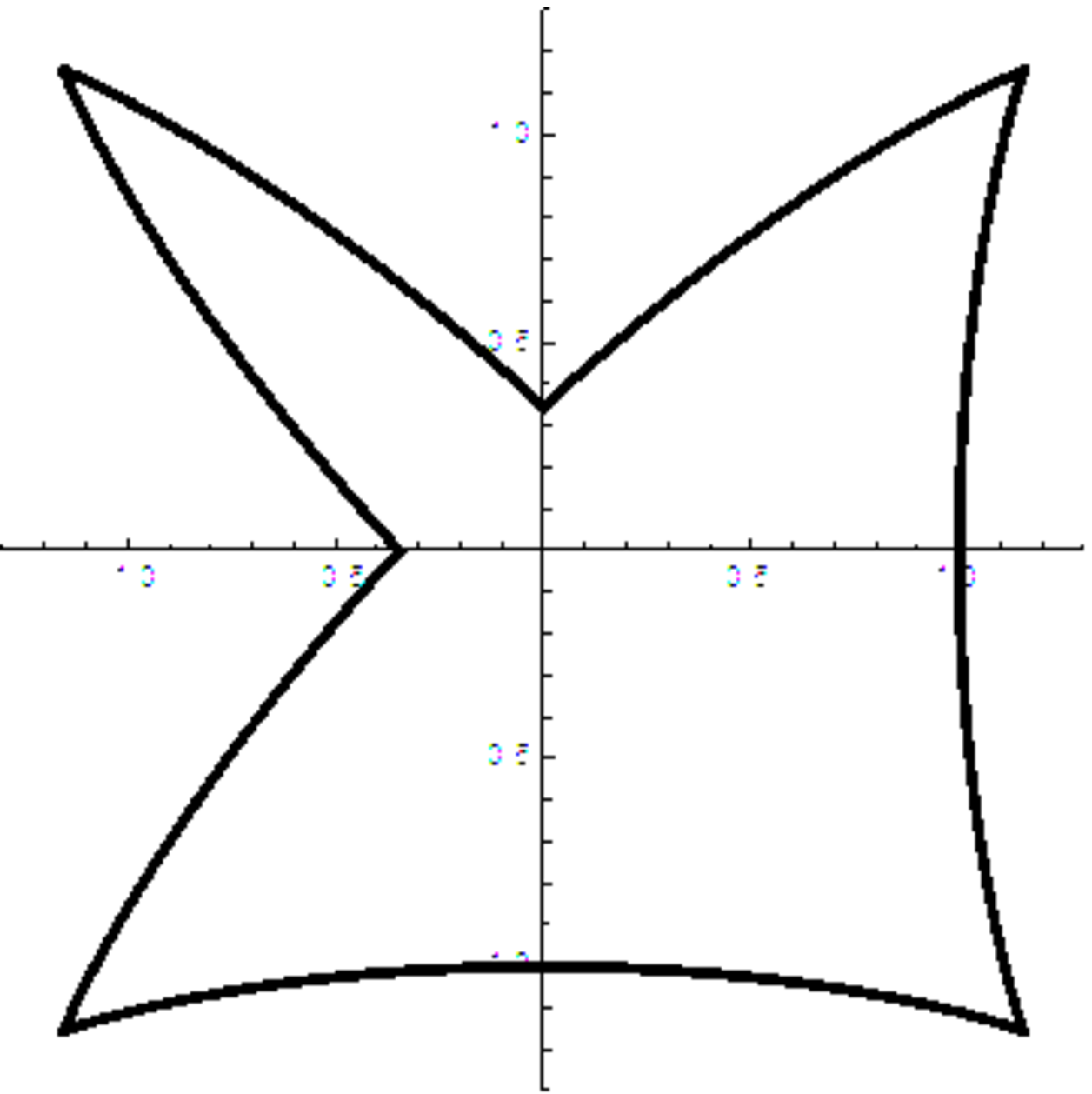} 
        \subcaption{CAMC curve}
\label{fig-2g}
  \end{minipage}   
  \hspace{5mm}     
    \begin{minipage}[b]{0.30\linewidth}
    \centering      \includegraphics[width=25mm,height=24mm,angle=0]{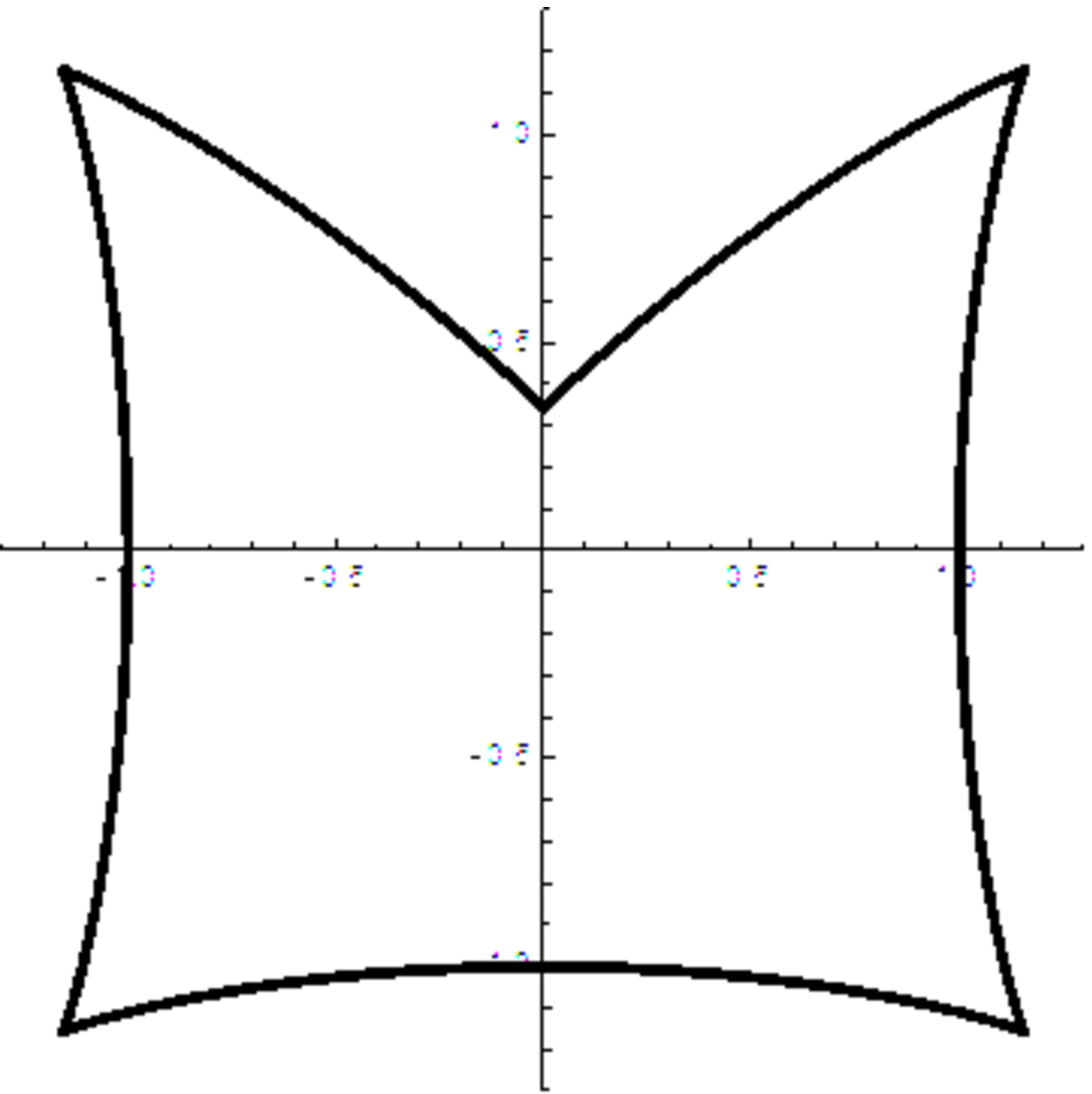} 
        \subcaption{CAMC curve}
\label{fig-2h}
  \end{minipage} 
  
  \caption{Closed CAMC curves  which are subsets of $\xi_\gamma(S^1)$ for $\gamma$ defined by (\ref{exk2017-1}) (cf. Figure \ref{fig:ex1-1}). For all of these eight curves, the anisotropic (mean) curvature for the outward-pointing normal is $-1$. They and the five curves in Figures \ref{fig-a}, \ref{fig-b}, \ref{fig-c}, \ref{fig-d}, \ref{fig-e} 
  give all closed CAMC curves included in $\xi_\gamma(S^1)$ (up to congruence in ${\mathbb R}^2$).
  }\label{fig-3}
\end{figure}

\begin{remark}\label{rem4}{\rm
In Proposition \ref{prop}, we gave six closed piecewise-$C^\infty$ curves which are subsets of $\xi_\gamma(S^1)$. Five of them were CAMC and the other one was not CAMC. We should remark that there are more piecewise-$C^\infty$ CAMC closed curves and non-CAMC closed curves. Figure \ref{fig-3} gives all of the other closed CAMC curves included in $\xi_\gamma(S^1)$ (up to congruence in ${\mathbb R}^2$).
}
\end{remark}

\section{Higher dimensional examples}\label{ex2}

In this section we give two higher dimensional examples by rotating $\gamma$ which was defined by (\ref{exk2017-1}) and studied in \S \ref{ex1} details. 

\begin{figure}[H]
 \centering
  \begin{minipage}[b]
  {0.20\linewidth}    
\centering \includegraphics[width=20mm,height=20mm,angle=0]{CHF180327SC.eps}
    \subcaption{The image of $\xi_{\gamma_1^0}$}\label{fig2A}
  \end{minipage}
  \hspace{5mm}
    \begin{minipage}[b]{0.25\linewidth}
    \centering
    \includegraphics[width=35mm,height=35mm,angle=0]{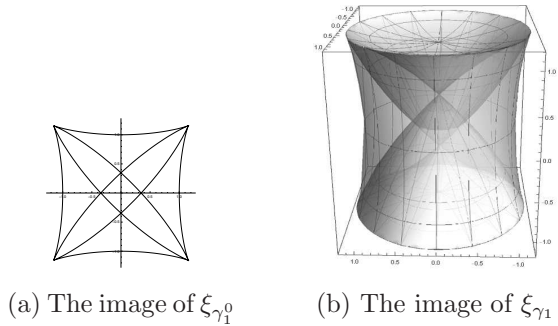}
    \subcaption{The image of $\xi_{\gamma_1}$}\label{fig3A}
  \end{minipage} 
  \caption{(a): The image of the Cahn-Hoffman map $\xi_{\gamma_1^0}:S^1\to {\mathbb R}^2$ for $\gamma_1^0:S^1 \to {\mathbb R}_{>0}$. (b): The image of the Cahn-Hoffman map $\xi_{\gamma_1}:S^2\to {\mathbb R}^3$ for $\gamma_1:S^2 \to {\mathbb R}_{>0}$.}\label{figA}
\end{figure}

Regarding (\ref{gmb0}), we consider the function  
$$
\overline\gamma_1^0(\nu_1, \nu_3)=\frac{\nu_1^6+\nu_3^6}
{(\nu_1^2+\nu_3^2)^{5/2}}
$$
defined on the $(\nu_1, \nu_3)$-plane. 
We denote the restriction of $\overline\gamma_1^0$ to $S^1$ by $\gamma_1^0$. 
Then, by Lemma \ref{lem1}, 
the Cahn-Hoffman map $\xi_{\gamma_1^0}:S^1 \to {\mathbb R}^2$ for $\gamma_1^0$ is represented as follows (Figure \ref{fig2A}).
\begin{eqnarray}
\xi_{\gamma_1^0}(\nu)
&=&
\bigl(
\nu_1(\nu_1^6+6\nu_1^4\nu_3^2-5\nu_3^6), 
\nu_3(-5\nu_1^6+6\nu_1^2\nu_3^4+\nu_3^6)
\bigr) \label{cf2A}\\
&=&
\bigl((\cos\theta)(\cos^6\theta+6\cos^4\theta\sin^2\theta-5\sin^6\theta), \nonumber\\
&&(\sin\theta)(-5\cos^6\theta+6\cos^4\theta\sin^2\theta+\sin^6\theta)
\bigr).\label{cf2A2}
\end{eqnarray}
($\nu=(\nu_1, \nu_3)=(\cos\theta, \sin\theta) \in S^1$).
The higher dimensional example obtained by rotating $\overline\gamma_1^0$ around the $\nu_3$-axis is given by
\begin{equation}\label{ex3A}
\overline\gamma_1(\nu_1, \nu_2, \nu_3)=\frac{(\nu_1^2+\nu_2^2)^3+\nu_3^6}
{(\nu_1^2+\nu_2^2+\nu_3^2)^{5/2}}, \quad (\nu_1, \nu_2, \nu_3) \in {\mathbb R}^3.
\end{equation}
The corresponding Cahn-Hoffman map $\xi_{\gamma_1}:S^2 \to {\mathbb R}^3$ is given as follows (Figure \ref{fig3A}). 
\begin{eqnarray}
\xi_{\gamma_1}(\nu)
&=&
\bigl((\cos\theta)(\cos^6\theta+6\cos^4\theta\sin^2\theta-5\sin^6\theta)(\cos\rho), \nonumber\\
&&(\cos\theta)(\cos^6\theta+6\cos^4\theta\sin^2\theta-5\sin^6\theta)(\sin\rho), \nonumber\\
&&(\sin\theta)(-5\cos^6\theta+6\cos^4\theta\sin^2\theta+\sin^6\theta)
\bigr),\label{cf3A}
\end{eqnarray}
($\nu=(\cos\theta\cos\rho, \cos\theta\sin\rho, \sin\theta) \in S^2$).

\begin{figure}[H]
 \centering
    \begin{minipage}[b]{0.25\linewidth}
    \centering   \includegraphics[width=15mm,height=15mm,angle=0]{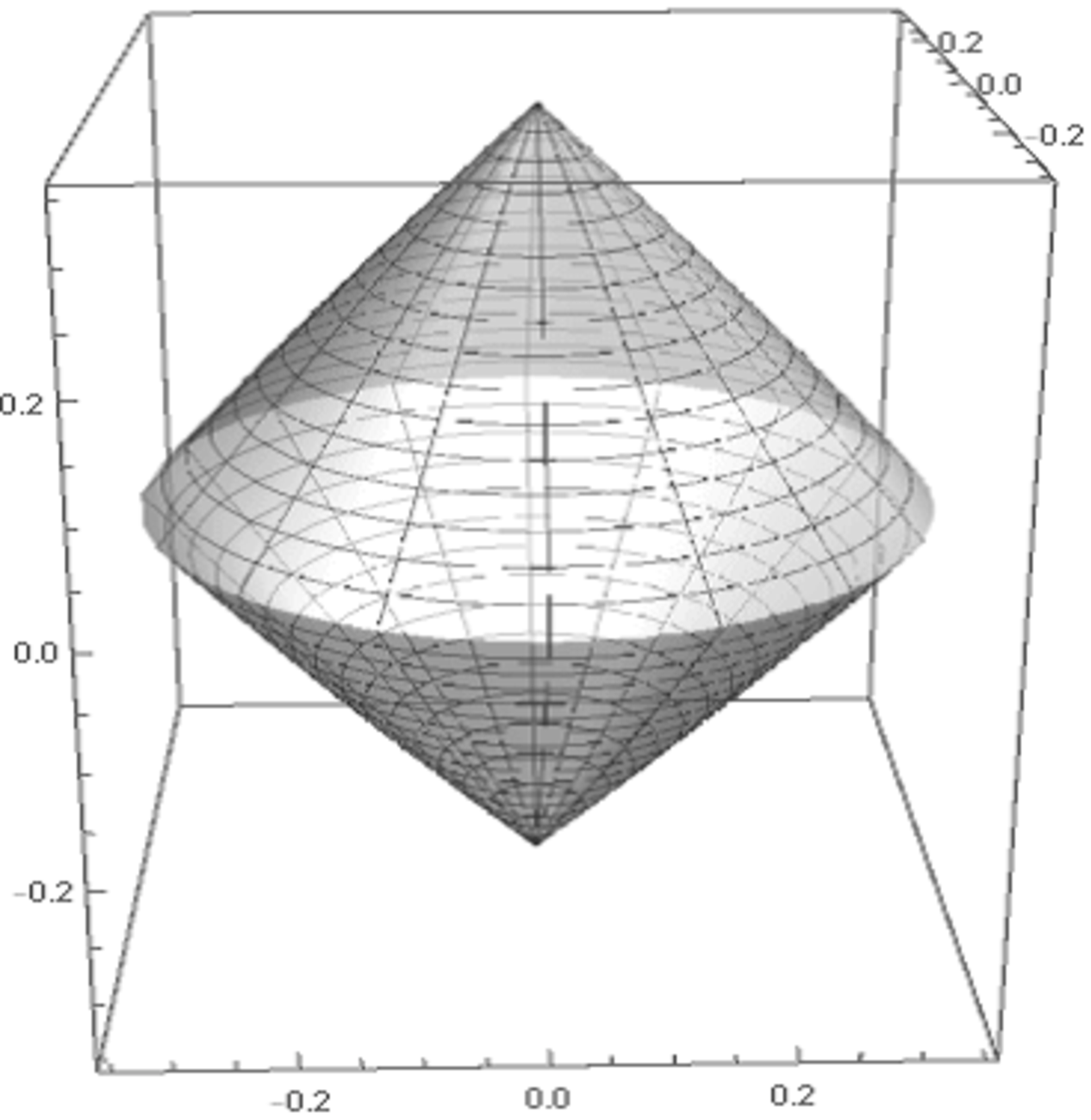}
    \subcaption{Wulff shape $W_{\gamma_1}$}\label{fig-3a}
  \end{minipage}
  \hspace{10mm} 
    \begin{minipage}[b]{0.25\linewidth}
    \centering    \includegraphics[width=25mm,height=25mm,angle=0]{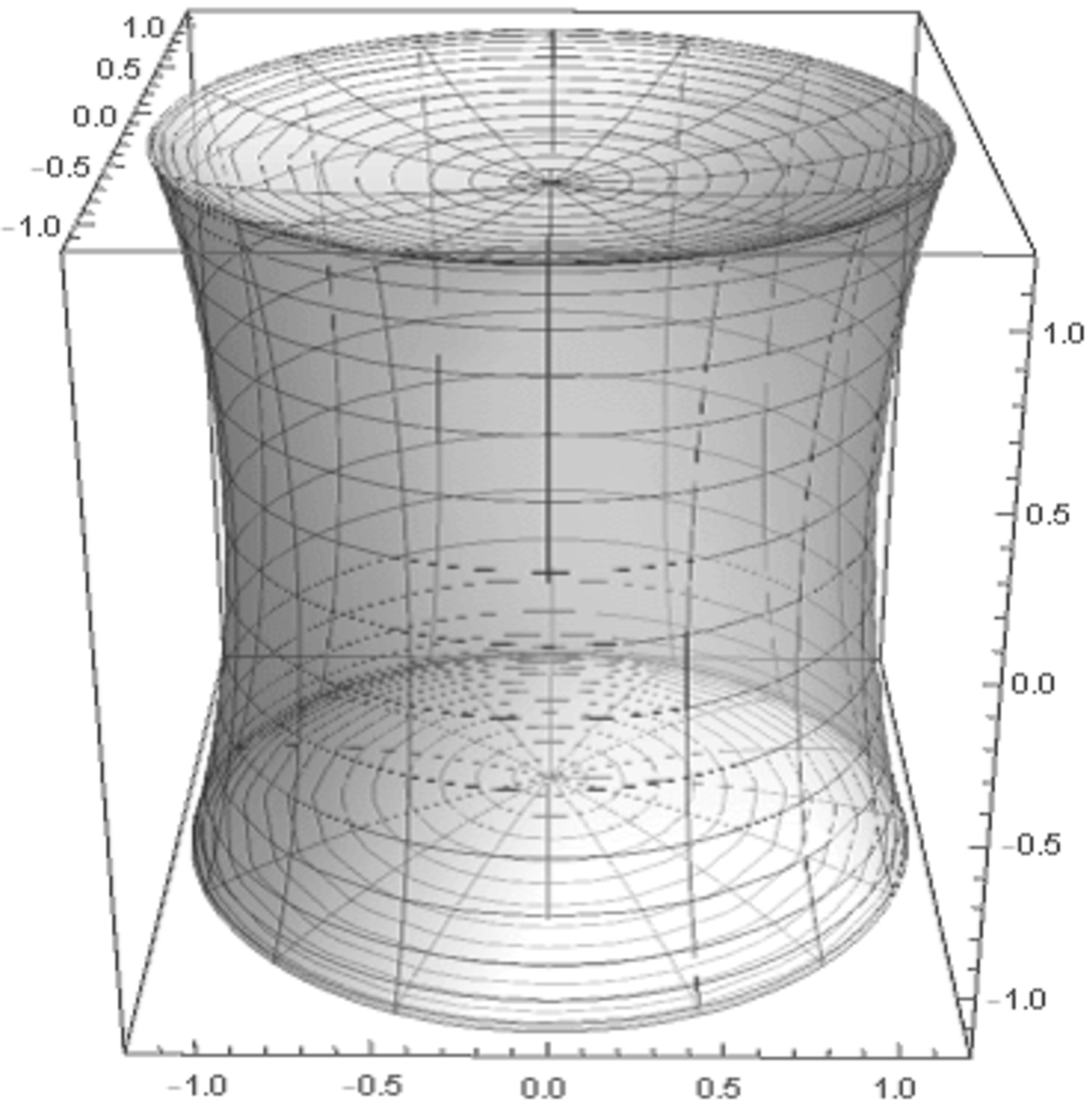}
    \subcaption{CAMC surface}\label{fig-3b}
  \end{minipage}   
\hspace{5mm}
\vskip1cm 
    \begin{minipage}[b]{0.25\linewidth}
    \centering    \includegraphics[width=27mm,height=26mm,angle=0]{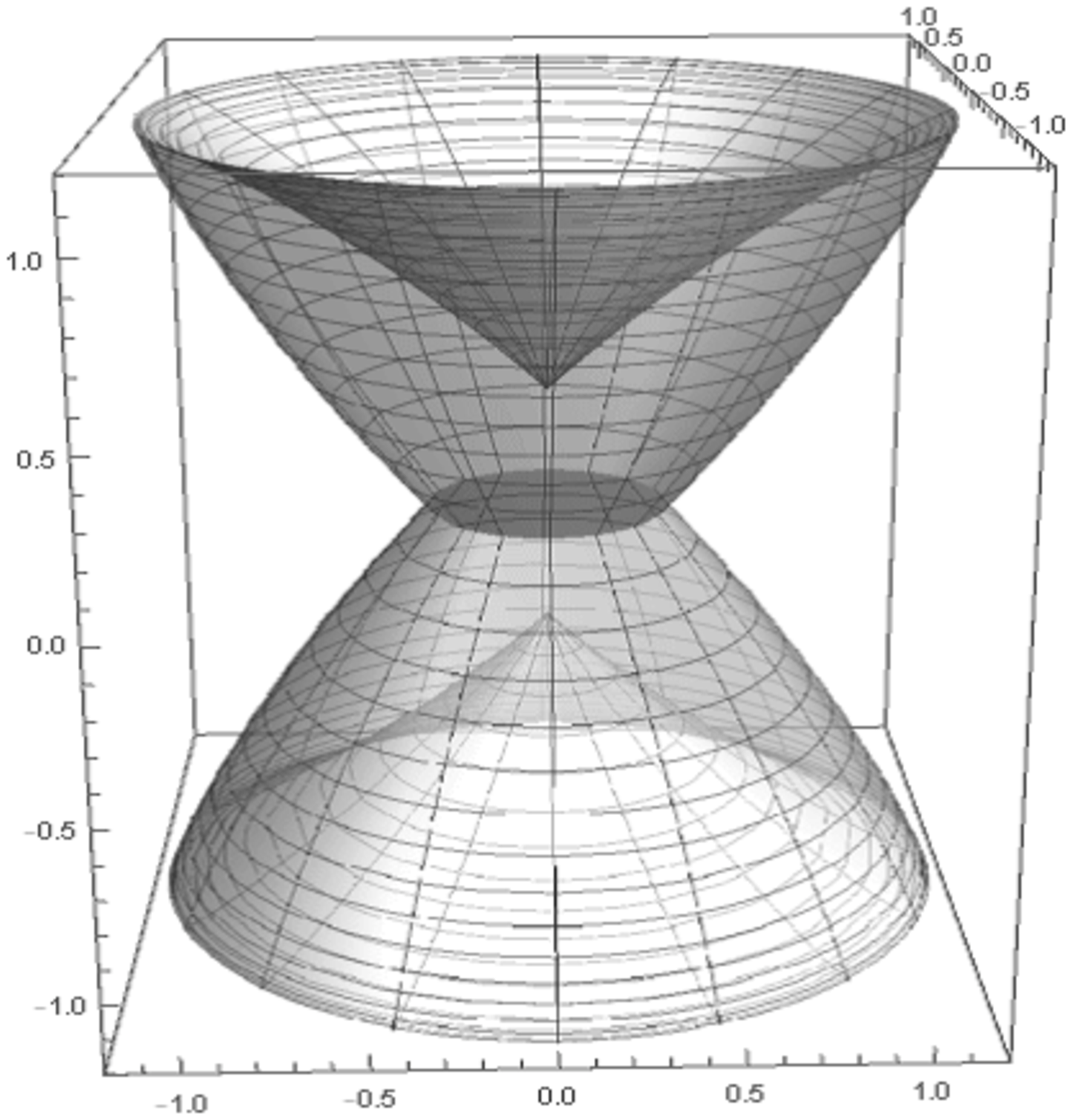}
    \subcaption{CAMC surface}\label{fig-3c}
  \end{minipage}   
\hspace{10mm} 
    \begin{minipage}[b]{0.30\linewidth}
    \centering      \includegraphics[width=31mm,height=26mm,angle=0]{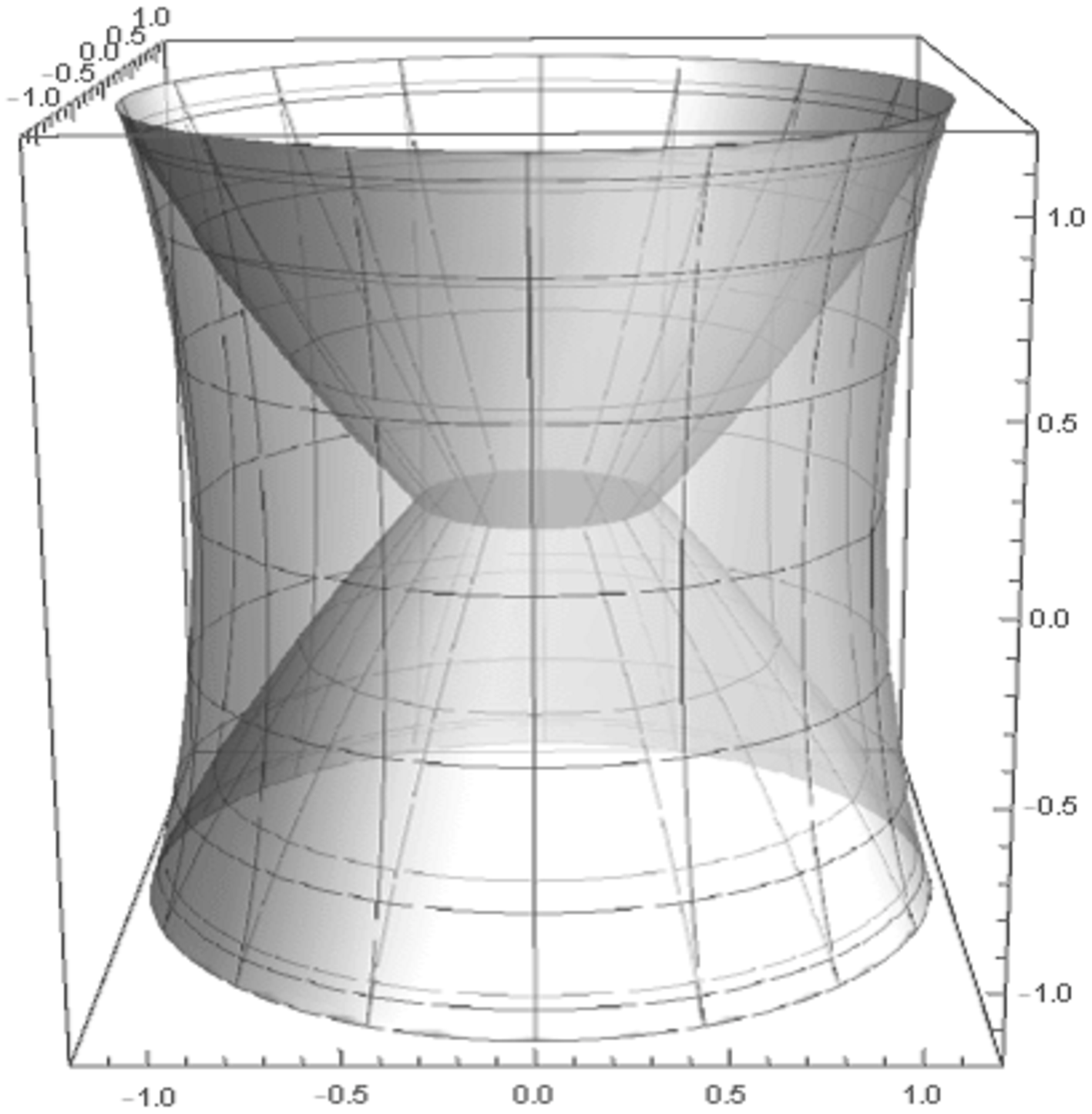} 
        \subcaption{not CAMC}\label{fig-3d}
  \end{minipage} 
  \caption{Some of the closed surfaces which are subsets of $\xi_{\gamma_1}(S^2)$ for $\gamma_1$ defined by (\ref{ex3AA}) (Figure \ref{fig3A}). They are surfaces given by rotating the curves $W_\gamma$, $(C_\gamma)_1$, $(C_\gamma)_4$, $(C_\gamma)_5$, respectively. (a), (b), (c): The anisotropic mean curvature for the outward-pointing normal is $-1$. (d) : The anisotropic mean curvature is $-1$ on the `outer part', while it is $1$ on the `inner part'. Hence, this surface is not CAMC.}\label{fig-q}
\end{figure}

By the same way as in \S \ref{ex1}, we get closed piecewise-$C^\infty$ CAMC surfaces and closed piecewise-$C^\infty$ non-CAMC surfaces for $\gamma_1$ which are subsets of $\xi_{\gamma_1}(S^2)$ and which are not the Wulff shape $W_{\gamma_1}$ (up to homothety and translation). 
In fact, we have the following:

\begin{proposition}\label{prop2}
Consider  
$\gamma_1:S^2\to{\mathbb R}_{>0}$ defined by 
\begin{equation}\label{ex3AA}
\gamma_1(\nu_1, \nu_2, \nu_3)=(\nu_1^2+\nu_2^2)^3+\nu_3^6, \quad (\nu_1, \nu_2, \nu_3) \in S^2.
\end{equation}
The Wulff shape $W_{\gamma_1}$  is the surface of revolution (Figure \ref{fig-3a}) given by
 rotating $W_{\gamma_1^0}$ (Figure \ref{fig-a}) around the vertical axis. 
 The two piecewise-$C^\infty$ closed surfaces (Figures \ref{fig-3b}, \ref{fig-3c}) given by rotating the closed curves $(C_\gamma)_1$, $(C_\gamma)_4$ (Figures \ref{fig-b}, \ref{fig-e}) around the vertical axis are CAMC. The piecewise-$C^\infty$ closed surface (Figures \ref{fig-3d}) given by rotating the closed curve $(C_\gamma)_5$ (Figure \ref{fig-f}) is not CAMC. 
\end{proposition}

\begin{figure}[H]
 \centering
  \begin{minipage}[b]
  {0.20\linewidth}    
\centering \includegraphics[width=28mm,height=28mm,angle=0]{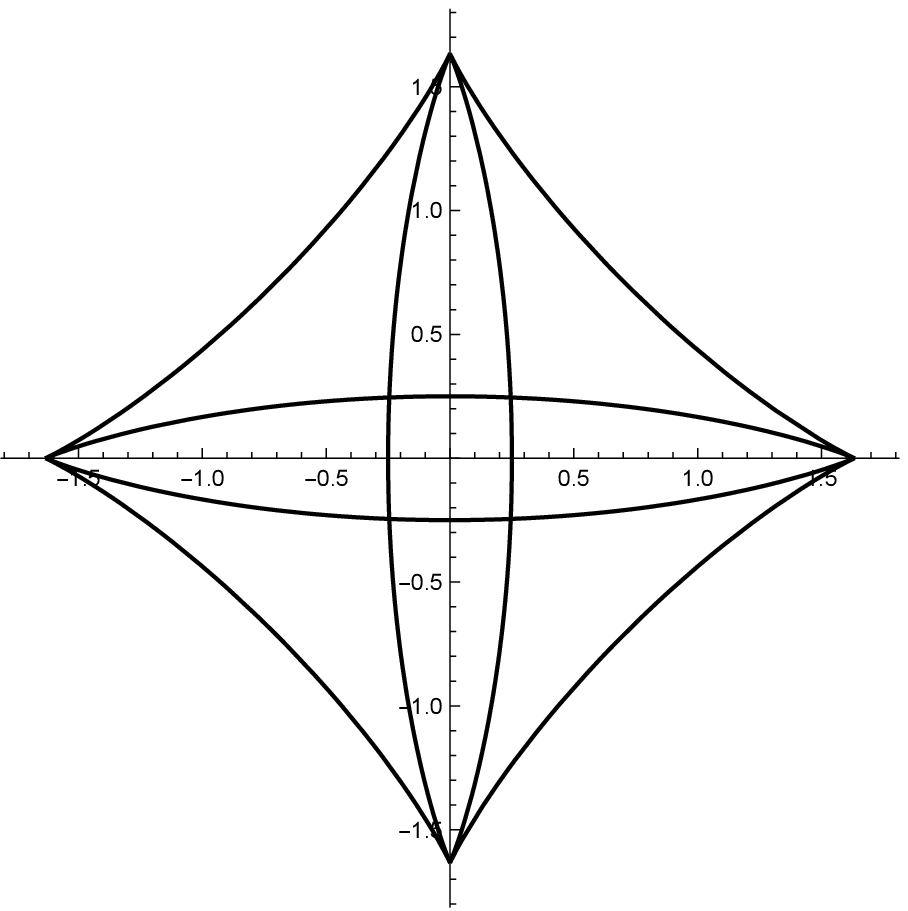}
    \subcaption{The image of $\xi_{\gamma_2^0}$}\label{fig2B}
  \end{minipage}
  \hspace{5mm}
    \begin{minipage}[b]{0.25\linewidth}
    \centering    \includegraphics[width=35mm,height=35mm,angle=0]{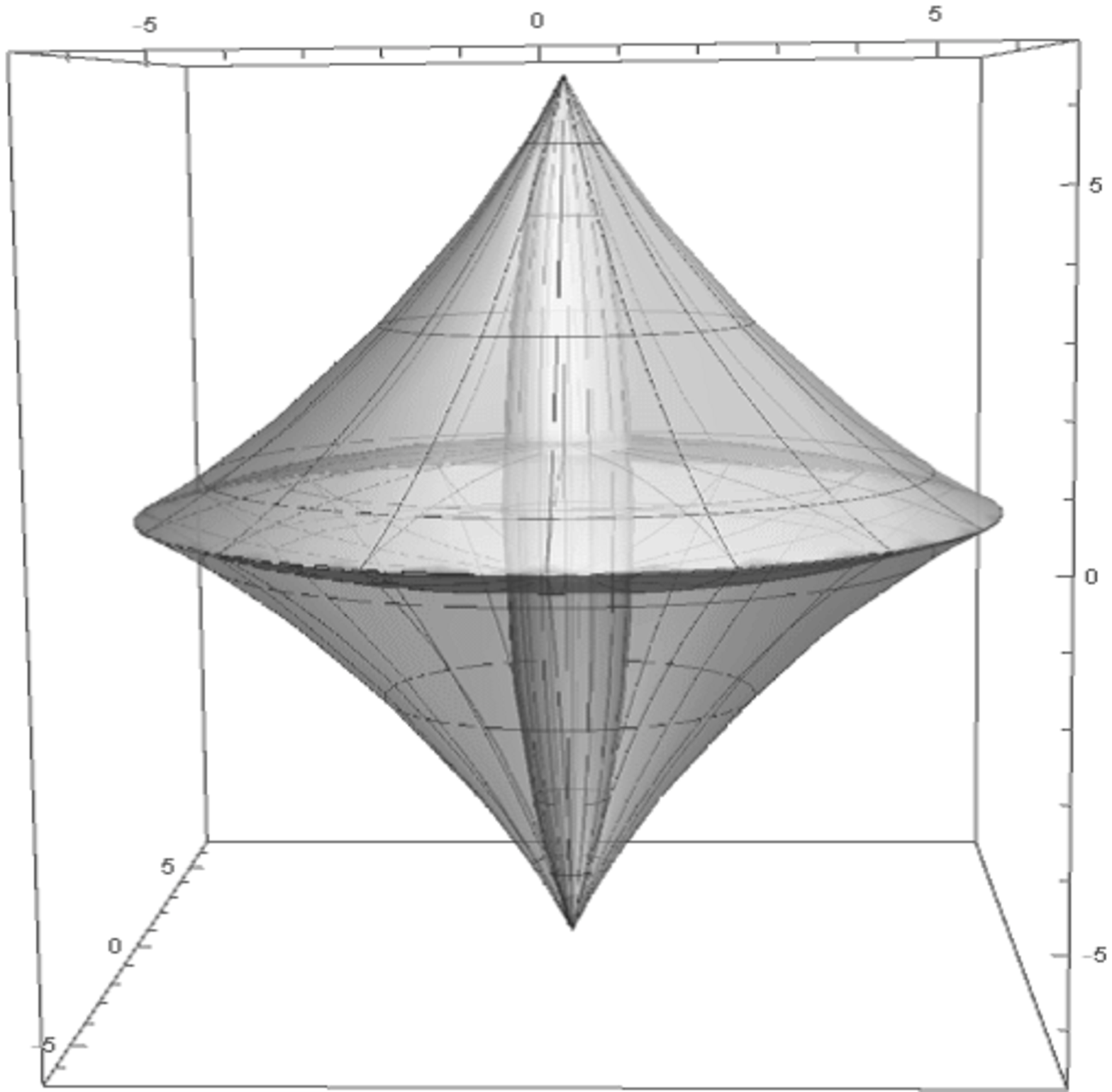}
    \subcaption{The image of $\xi_{\gamma_2}$}\label{fig3B}  
      \end{minipage} 
      \caption{(a): The image of the Cahn-Hoffman map $\xi_{\gamma_2^0}:S^1\to {\mathbb R}^2$ for $\gamma_2^0:S^1 \to {\mathbb R}_{>0}$. (b): The image of the Cahn-Hoffman map $\xi_{\gamma_2}:S^2\to {\mathbb R}^3$ for $\gamma_2:S^2 \to {\mathbb R}_{>0}$.}\label{figB}
\end{figure}

Let us give another example. 
We rotate $\overline\gamma_1^0$ around the origin by $\pi/4$ and obtain the following example. 
\begin{equation}\label{eq22}
\overline\gamma_2^0(\nu_1, \nu_3)=\frac{\nu_1^6+15\nu_1^4\nu_3^2+15\nu_1^2\nu_3^4+\nu_3^6}
{4(\nu_1^2+\nu_3^2)^{5/2}}
\end{equation}
defined on the $(\nu_1, \nu_3)$-plane. 
We denote the restriction of $\overline\gamma_2^0$ to $S^1$ by $\gamma_2^0$. 
The Cahn-Hoffman map $\xi_{\gamma_2^0}:S^1 \to {\mathbb R}^2$ for $\gamma_2^0$ is obtained by rotating $\xi_{\gamma_1^0}$ around the origin by $\pi/4$, and so it 
is represented as follows (Figure \ref{fig2B}). 
\begin{eqnarray}
\xi_{\gamma_2^0}(\nu)
\!\!\!&=&\!\!\!
\frac{1}{\;4\;}\bigl(
\nu_1(\nu_1^6-9\nu_1^4\nu_3^2+15\nu_1^2\nu_3^4+25\nu_3^6), 
\nu_3(25\nu_1^6+15\nu_1^4\nu_3^2-9\nu_1^2\nu_3^4+\nu_3^6)
\bigr) \label{cf2B}\\
\!\!\!&=&\!\!\!
\frac{1}{\;4\;}\bigl((\cos\theta)(\cos^6\theta-9\cos^4\theta\sin^2\theta+15\cos^2\theta\sin^4\theta+25\sin^6\theta), \nonumber\\
&& \ \ \ \ (\sin\theta)(25\cos^6\theta+15\cos^4\theta\sin^2\theta-9\cos^2\theta\sin^4\theta+\sin^6\theta)
\bigr).\label{cf2B2}
\end{eqnarray}
($\nu=(\nu_1, \nu_3)=(\cos\theta, \sin\theta) \in S^1$).
The higher dimensional example obtained by rotating $\overline\gamma_2^0$ around the $\nu_3$-axis is given by
\begin{equation}\label{ex3B}
\overline\gamma_2(\nu_1, \nu_2, \nu_3)=\frac{(\nu_1^2+\nu_2^2)^3+15(\nu_1^2+\nu_2^2)^2\nu_3^2+15(\nu_1^2+\nu_2^2)\nu_3^4+\nu_3^6}
{4(\nu_1^2+\nu_2^2+\nu_3^2)^{5/2}}, \quad (\nu_1, \nu_2, \nu_3) \in {\mathbb R}^3.
\end{equation}
The restriction $\gamma_2$ of $\overline\gamma_2:{\mathbb R}^3 \to {\mathbb R}_{\ge 0}$ to $S^2$ can be written as
\begin{equation}\label{ex3BB}
\gamma_2(\nu_1, \nu_2, \nu_3)=(\nu_1^2+\nu_2^2)^3+15(\nu_1^2+\nu_2^2)^2\nu_3^2+15(\nu_1^2+\nu_2^2)\nu_3^4+\nu_3^6, \quad (\nu_1, \nu_2, \nu_3) \in S^2.
\end{equation}
The corresponding Cahn-Hoffman map $\xi_{\gamma_2}:S^2 \to {\mathbb R}^3$ is given as follows (Figure \ref{fig3B}). 
\begin{eqnarray}
\xi_{\gamma_2}(\nu)
\!\!\!&=&\!\!\!
\frac{1}{\;4\;}\bigl((\cos\theta)(\cos^6\theta-9\cos^4\theta\sin^2\theta+15\cos^2\theta\sin^4\theta+25\sin^6\theta)(\cos\rho), \nonumber\\
&& \ \ \ \ (\cos\theta)(\cos^6\theta-9\cos^4\theta\sin^2\theta+15\cos^2\theta\sin^4\theta+25\sin^6\theta)(\sin\rho), \nonumber\\
&&\ \ \ \ (\sin\theta)(25\cos^6\theta+15\cos^4\theta\sin^2\theta-9\cos^2\theta\sin^4\theta+\sin^6\theta)
\bigr),\label{cf3B}
\end{eqnarray}
($\nu=(\cos\theta\cos\rho, \cos\theta\sin\rho, \sin\theta) \in S^2$).

\begin{figure}[H]
 \centering
   \begin{minipage}[b]
  {0.20\linewidth}    
\centering \includegraphics[width=14mm,height=11mm,angle=0]{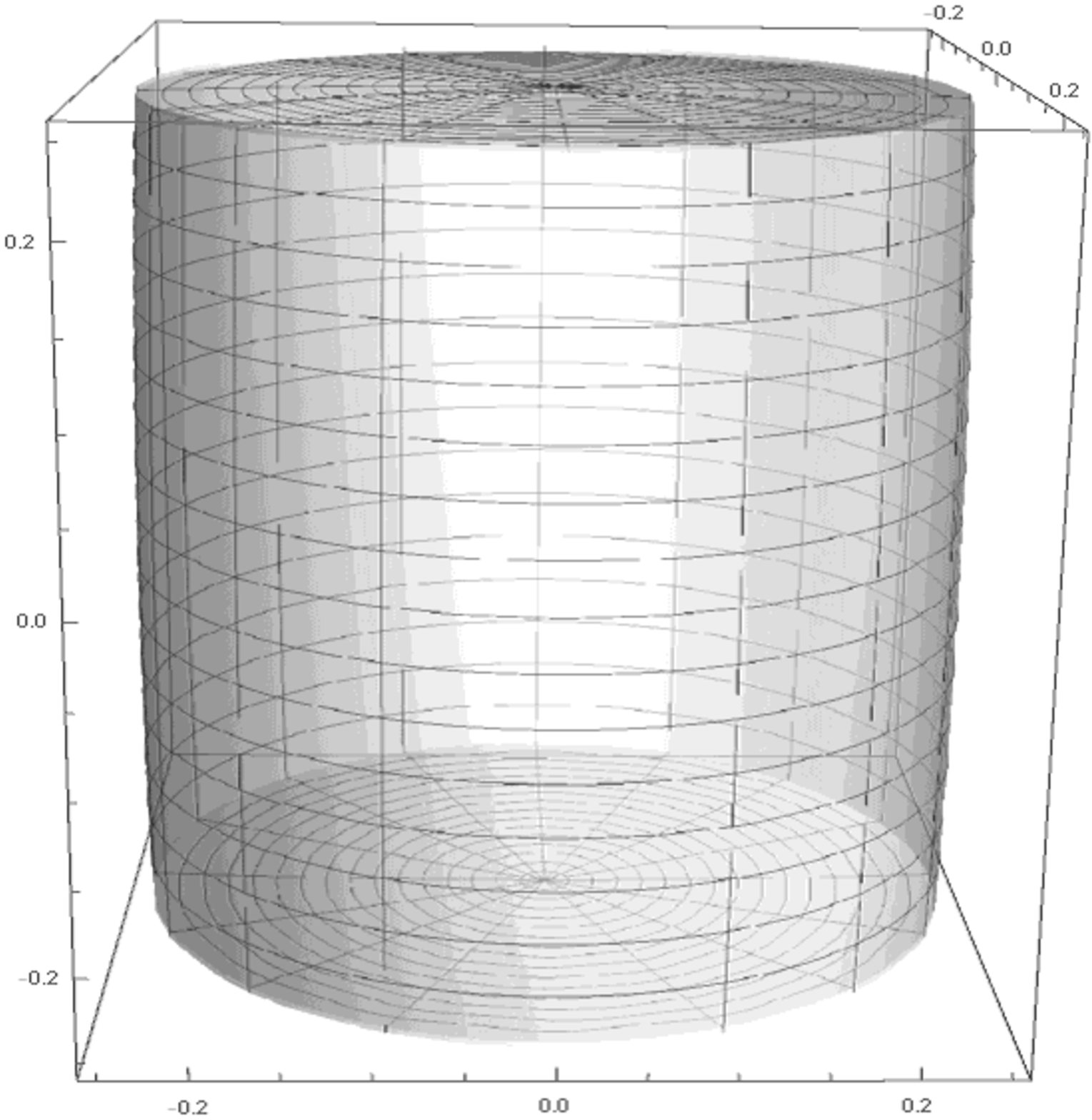}
    \subcaption{Wulff shape $W_{\gamma_2}$}\label{fig1C}
  \end{minipage}
  \begin{minipage}[b]
  {0.50\linewidth}    
\centering \includegraphics[width=33mm,height=27mm,angle=0]{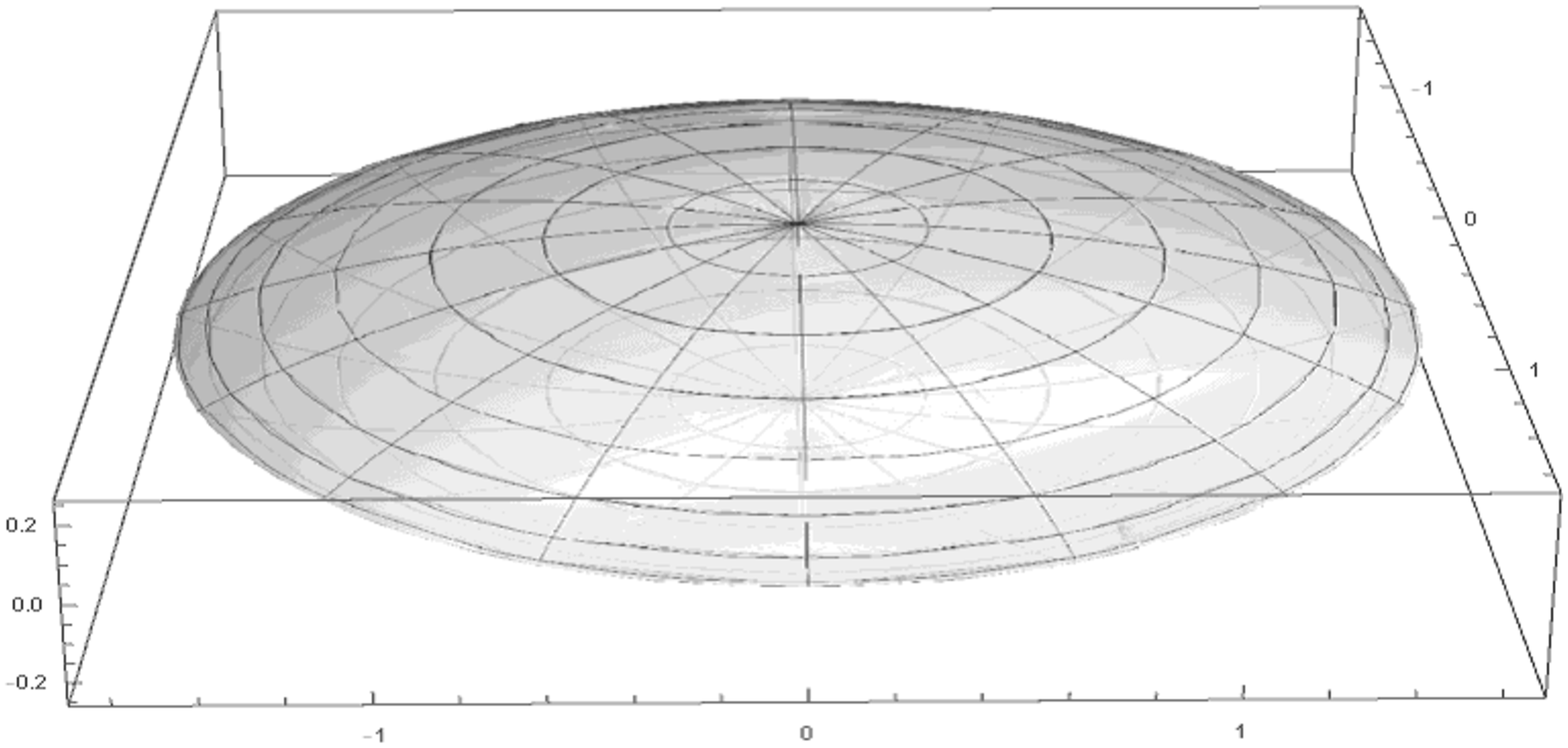}
    \subcaption{CAMC surface for $\gamma_2$}\label{fig2C}
  \end{minipage}
    \begin{minipage}[b]{0.25\linewidth}
    \centering    \includegraphics[width=33mm,height=27mm,angle=0]{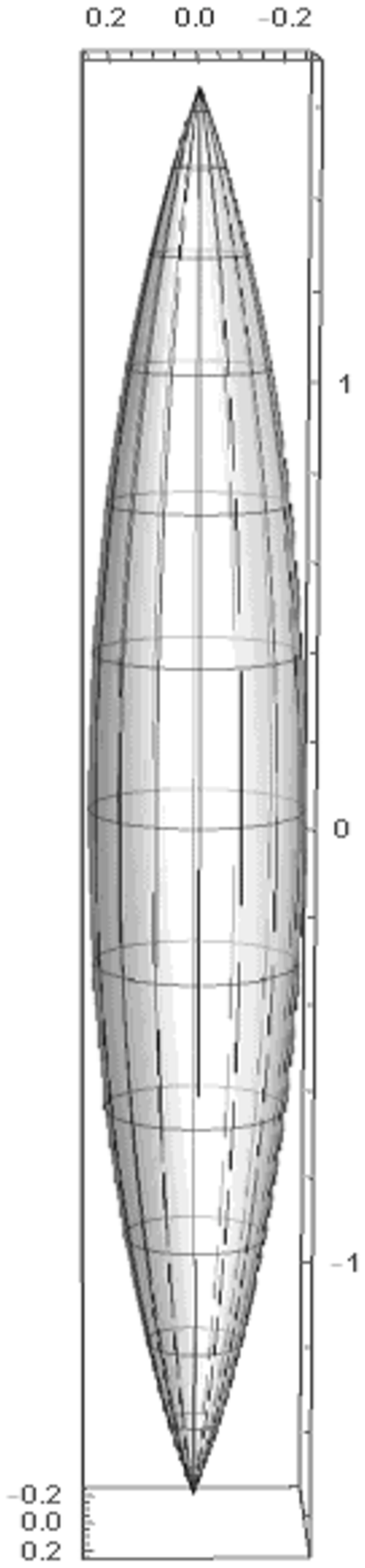}
    \subcaption{CAMC surface for $\gamma_2$}\label{fig3C}  
      \end{minipage} 
        \caption{Some of the closed surfaces which are subsets of $\xi_{\gamma_2}(S^2)$ for $\gamma_2$ defined by (\ref{ex3BB}) (Figure \ref{fig3B}). The anisotropic mean curvature for the outward-pointing normal is $-1$. }\label{figC}
\end{figure}

By the same way as above, we get closed piecewise-$C^\infty$ CAMC surfaces (Figure \ref{figC}) and closed piecewise-$C^\infty$ non-CAMC surfaces for $\gamma_2$ which are subsets of $\xi_{\gamma_2}(S^2)$ and which are not the Wulff shape $W_{\gamma_2}$ (up to homothety and translation).

\section{Proofs of Theorems \ref{nonuni1}, \ref{nonuni2}}\label{pf}

The functions $\gamma_1^0:S^1\to{\mathbb R}_{>0}$, $\gamma_2^0:S^1\to{\mathbb R}_{>0}$ and their rotations around the vertical axis give examples which prove Theorems \ref{nonuni1}, \ref{nonuni2}. 
In fact, Propositions \ref{prop}, \ref{prop2} give suitable examples for $n=1, 2$, respectively. Also, higher dimensional examples are obtained by the method given in \S \ref{ex2}.

\section{Applications to anisotropic mean curvature flow:  proofs of Theorems \ref{ASS},  \ref{AMCF}} \label{app}

Let $\gamma:S^n\to{\mathbb R}_{>0}$ be of $C^2$ with Cahn-Hoffman map $\xi_\gamma$. 
Let $X_t:M \to {\mathbb R}^{n+1}$ be one-parameter family of embedded piecewise-$C^2$ hypersurfaces with anisotropic mean curvature $\Lambda_t$. Assume that the Cahn-Hoffman field $\tilde{\xi}_t$ along $X_t$ is defined on $M$. 
If $X_t$ satisfies 
$\partial X_t/\partial t = \Lambda_t \tilde{\xi}_t$, 
it is called an anisotropic mean curvature flow, which diminishes the anisotropic energy if $\Lambda_t \not\equiv0$. In fact,
\begin{eqnarray}\label{ft2}
\frac{d{\cal F}_\gamma(X_t)}{dt}
\!\!&=&\!\! -\int_M n\Lambda_t \Big\langle 
\frac{\partial X_t}{\partial t}, \nu_t
\Big\rangle \;dA_t
\!=\! -\int_M n\Lambda_t^2 \langle 
D\gamma+\gamma(\nu_t)\nu_t, \nu_t
\rangle \;dA_t \nonumber \\
&=& -\int_M n\Lambda_t^2\gamma(\nu_t)\;dA_t
\le 0
\end{eqnarray}
holds. 

\vskip0.5truecm

\noindent{\it Proof of Theorem \ref{ASS}.} \ 
Since the anisotropic mean curvature of $\xi_\gamma$ is $-1$, $\displaystyle \Lambda_t=\frac{-1}{\sqrt{2(c-t)}}$ holds. On the other hand, $\tilde{\xi}_t=\xi_\gamma$ holds. These two facts imply that (i) and (ii) hold. 
\hfill $\Box$

\vskip0.5truecm

\noindent{\it Proof of Theorem \ref{AMCF}.} \ 
Examples stated in Propositions \ref{prop}, \ref{prop2} give the desired result. 
\hfill $\Box$



\begin{flushleft}
Yoshiki J{\footnotesize IKUMARU} \\
Graduate School of Mathematics\\
Kyushu University  \\
Motooka Nishi-ku, Fukuoka \\
FUKUOKA 819-0395,
JAPAN \\
E-mail: {\tt y-jikumaru@math.kyushu-u.ac.jp} \\
 \end{flushleft}
 
  \begin{flushleft}
Miyuki K{\footnotesize OISO} \\
Institute of Mathematics for Industry  \\
Kyushu University  \\
Motooka Nishi-ku, Fukuoka \\
FUKUOKA 819-0395,
JAPAN \\
E-mail: {\tt koiso@math.kyushu-u.ac.jp} \\
 \end{flushleft}

\end{document}